\documentclass[a4paper,11pt]{article}
\usepackage[english]{babel}
\usepackage[margin=2.8cm]{geometry}
\usepackage[noblocks]{authblk}

\usepackage[T1]{fontenc}
\usepackage[latin1]{inputenc}

\usepackage{amsmath,amsfonts,amssymb,amsthm,mathrsfs,enumerate,color}
\usepackage{graphicx}

\usepackage[draft]{fixme}
\fxusetheme{color}
\FXRegisterAuthor{b}{ab}{B}
\FXRegisterAuthor{pe}{ape}{Pe}
\FXRegisterAuthor{p}{ap}{P}
\FXRegisterAuthor{s}{as}{S}
\usepackage{tikz}

\newcommand{\cP}{\mathcal{P}}
 
 \newcommand{\PP}{\mathbb{P}}
 \newcommand{\Nset}{\mathbb{N}}
 \newcommand{\eps}{\varepsilon}
 \newcommand{\lpref}{\smash{\raisebox{3.5pt}{\!\!\!\begin{tabular}{c}$\hskip-4pt\scriptstyle\longleftarrow$ \\[-7pt]{\rm pref}\end{tabular}\!\!}}}
\newcommand{\petitlpref}{\smash{\raisebox{2.5pt}{\!\!\!\begin{tabular}{c}$\hskip-2pt\scriptscriptstyle\longleftarrow$ \\[-9pt]{$\scriptstyle\hskip 1pt\rm pref$}\end{tabular}\!\!}}}
\newcommand\1{\leavevmode\hbox{\rm \small1\kern-0.35em\normalsize1}}
\newcommand\ind[1]{\1_{\left\{#1\right\}}}
\def\ie{\emph{i.e.~}}
\def\g#1{\mathbb #1}
\def\rond#1{\mathcal #1}
\newcommand{\pff}{\noindent {\sc Proof.}\ }
\def\QED{\hfill\vrule height 1.5ex width 1.4ex depth -.1ex \vskip 10pt}

 \newtheorem{thm}{Theorem}[section]
 \newtheorem{cor}[thm]{Corollary}
 \newtheorem{lem}[thm]{Lemma}
 \newtheorem{prop}[thm]{Proposition}
 \newtheorem{rem}[thm]{Remark}
 \newtheorem{defi}[thm]{Definition}
\title{Persistent random walks, variable length Markov chains and piecewise deterministic Markov processes %\footnote{
\footnote{
Supported in part by the  ANR grants :  "Malliavin, Stein and Stochastic Equations with Irregular
Coefficients" [ANR-10-BLAN-0121], "Mathematical Analysis of Neuronal Dynamics"  [ANR-09-BLAN-0008-01] and "Approches Spatio­-Temporelles pour la modélisation du Risque" [ANR- 08-BLAN-0314-01].
}}
 \author[1]{Peggy {\sc Cénac}}
 \author[2]{Brigitte {\sc Chauvin}}
 \author[1]{Samuel {\sc Herrmann}}
 \author[3]{Pierre {\sc Vallois}}
 \affil[1]{Institut de Math\'ematiques de Bourgogne (IMB) - UMR CNRS 5584,

 Universit\'e de Bourgogne,  B.P. 47 870

   21078 Dijon Cedex, France}
 \affil[2]{Laboratoire de Math\'ematiques de Versailles (LMV) - UMR CNRS 8100
 
45 avenue des Etats-Unis

78035 Versailles Cedex, France}
 \affil[3]{Université de Lorraine, Institut de Mathématiques Elie Cartan, 
 
 INRIA-BIGS, CNRS UMR 7502,
 
BP 239, F-54506 Vandoeuvre-lès-Nancy Cedex, France}
\begin{document}
%La version de TikZ est : \pgfversion
\maketitle
\begin{abstract}
A classical random walk $(S_t,\, t\in\mathbb{N})$ is defined by $S_t:=\displaystyle\sum_{n=0}^t X_n$, where $(X_n)$ are i.i.d.
When the increments $(X_n)_{n\in\mathbb{N}}$ are a one-order Markov chain,
a short memory is introduced in the dynamics of $(S_t)$.
This so-called ``persistent'' random walk is nolonger Markovian and, under suitable conditions, the rescaled
process converges towards the integrated telegraph noise (ITN) as the time-scale
and space-scale parameters tend to zero (see \cite{Herrmann-Vallois, Tapiero-Vallois, Tapiero-Vallois2}). The ITN process is effectively non-Markovian too. 
The aim is to consider persistent random walks $(S_t)$ whose increments are Markov chains with variable order which can be infinite. This variable memory is enlighted by a one-to-one correspondence between $(X_n)$ and a suitable Variable Length Markov Chain (VLMC), since for a VLMC the dependency from the past can be unbounded. 
 The key fact is to consider the non Markovian letter process $(X_n)$ as the margin of a couple $(X_n,M_n)_{n\ge 0}$ where $(M_n)_{n\ge 0}$ stands for the memory of the process $(X_n)$. We prove that, under a suitable rescaling, $(S_n,X_n,M_n)$ converges in distribution towards a time continuous process $(S^0(t),X(t),M(t))$. The process $(S^0(t))$ is a semi-Markov and Piecewise Deterministic Markov Process whose paths are piecewise linear.
\end{abstract}

Variable length Markov chain.
Integrated telegraph noise.
Piecewise Deterministic Markov Processes. 
Semi Markov processes. 
Variable memory. 
Simple and double infinite combs.
\section{Introduction}
\label{intro}
Classical random walks are defined by 
\begin{equation}
\label{def-persist-part}
S_t:=\displaystyle\sum_{n=0}^t X_n, 
\end{equation}
for $t\in\g N$ and for i.i.d.\@ increments $(X_n)_{n\in\mathbb{N}}$. It is well known that a suitable rescaling of the random walk permits to obtain the standard Brownian motion as the time-scale and space-scale parameters tend to zero.  When the increments $(X_n)_{n\in\mathbb{N}}$ are defined as a one-order Markov chain, a short memory in the dynamics of the stochastic paths is introduced: the process is called in the literature the \emph{persistent} random walk or a correlated random walk or also a Kac walk (see \cite{eckstein00, renshaw81, weissbook94, weiss02}). The random walk is nolonger Markovian and, under suitable conditions, the rescaled
process converges towards the integrated telegraph noise (ITN), see  \cite{Herrmann-Vallois,Tapiero-Vallois} and \cite{Tapiero-Vallois2}. The ITN process is
effectively non-Markovian too.

Our aim is to define processes $(X_n)_{n\in\mathbb{N}}$ with variable memory and thus to generalize this
convergence result to random walks whose increments are higher
order Markov chains. When $(X_n)_{n\in\mathbb{N}}$ is a Markov chain of
finite order, it is natural to think that the limit process should be very close to the
integrated telegraph noise. That is why we are mostly interested in constructing
infinite length Markov chain or in dealing with \emph{Variable
Length Markov Chains} (VLMC) for which the dependency from the past is non bounded. 

A VLMC can be defined as follows (this probabilistic presentation comes from \cite{ccpp}, other more statistic points of view can be found in \cite{Rissanen, Galves-Locherbach}). 
Let $\mathcal{L}= \{0,1\}^{-\mathbb{N}}$ be the set of left-infinite words on the alphabet $\{0,1\}$. Consider a complete (each node has $0$ or $2$ children) binary tree whose finite leaves $\mathcal{C}$ are words on the alphabet $\{0,1\}$. To each leaf $c$ (not necessarily finite) is attached a Bernoulli distribution denoted by $q_{c}$. Each leaf is called a \emph{context} and this probabilized tree is called a \emph{context tree}. See for instance the simple infinite comb in Figure~\ref{Fig.arbre}: the set of leaves $\mathcal{C}$ is defined by
\[
\mathcal{C}:=\{0^n1,\, n\ge 0\}\cup \{0^\infty\}
\]
where $0^n1$ represents the sequence $00\ldots01$ composed with
$n$ characters $'0'$ and one character $'1'$. By convention $0^01=1$. The set of leaves
contains one infinite leaf $0^\infty$ and a countable set of
finite leaves $0^n1$. The \emph{prefix function} $\lpref :\mathcal{L}=\{0,1\}^{-\mathbb{N}}\to\mathcal{C}$ indicates the length of the last run of $'0'$: for instance,
 \[
 \lpref (\ldots 1000)=0001=0^31.
 \]
For a general context tree and for any left-infinite word $U$, we define $\lpref(U)$ in a similar way as the first suffix of $U$ reading from right to left appearing as a leaf of the context tree. The associated VLMC is the $\mathcal{L}$-valued Markov chain $(U_n)_{n\geqslant 0}$ defined by the transitions
\begin{equation}
 \label{eq:def:VLMC}
 \PP(U_{n+1} = U_n\ell| U_n)= q_{\petitlpref (U_n)}(\ell)
 \end{equation}
where $\ell\in \{0,1\}$ is any letter. Notice that the VLMC is entirely determined by the data $q_{c}, c \in \mathcal{C}$. Moreover the order of dependence  (the memory) depends on the past itself.

For a given VLMC $(U_n)_{n\geqslant 0}$, define $X_n$ as the last letter of $U_n$ for any $n\geq 0$. When the context tree associated with $(U_n)$ is infinite, then the letter process $(X_n)_{n\geqslant 0}$ is non Markovian, because the transition probabilities (\ref{eq:def:VLMC}) indicate that $X_{n+1}$ depends on a variable and unbounded number of previous letters. The corresponding random walk $(S_t)$ defined by \eqref{def-persist-part} is non Markovian anymore, it is somehow \emph{very persistent}, so we investigate the following natural questions:  
is the random walk of the same nature as in the one-order Markov case? Is the rescaled process convergent to some limit process? Is the limiting process analog to the ITN?

Recall that $X_n$ is the last letter of a VLMC $(U_n)$. The key point of view is the following: we consider the non Markovian letter process $(X_n)$ as the margin of a couple $(X_n,M_n)_{n\ge 0}$ where $(M_n)_{n\ge 0}$ stands for the memory of the process $(X_n)$.  It is reasonable to believe that $M_n=|\lpref (U_n)|$ is a good candidate, where the notation $|w|$ stands for the length of a word $w$. More precisely in the particular case of a two-letter alphabet $\mathcal{A}$, the Markov chain $(X_n, M_n)_{n\ge 0}$ valued in the state space $\mathcal{A}\times \mathbb{N}^*$ is defined by the transition probabilities: let $\ell, \ell '\in\rond A, \ell\not= \ell '$,
\begin{align*}
\left\{\begin{array}{l}
Q\Big((\ell,n),(\ell,n+1)\Big)=1-\alpha_{\ell,n},\\
Q\Big((\ell,n),(\ell ',1)\Big)=\alpha_{\ell,n}.
\end{array}\right.
\end{align*}
Note that $\alpha_{\ell,k}$ is the probability of changing letter after a run of length $k$ of letter $\ell$, that is
\begin{equation}\label{desc-alpha}
\alpha_{\ell,k}=\mathbb{P}(X_{n+1}\neq \ell\vert X_n=\ell,\ M_n=k).
\end{equation}
Introducing the sequence of breaking times:
\[
T_0=0,\quad T_{k+1}=\inf\{ n>T_k,\ X_n\neq X_{T_k} \}
\]
it is easy to see that $(X_n,T_n)_{n\ge 0}$ is a semi-Markov process (see \cite[Chapter 10]{cinlar} and \cite{Harlamov,janssen2006}).
In Section~\ref{sect:model}, we consider a Markov chain
$(X_n,M_n)_{n\ge 0}$, where $(X_n)_{n\geq 0}$ is a letter process, the letters belong to an alphabet $\rond A := \{a_1,a_2,\ldots,a_K\} $, and $(M_n)_{n\ge 0}$ stands for the memory of the process $(X_n)$.  The state space associated with this Markov chain is $\{a_1,a_2,\ldots,a_K\}\times\overline{\mathbb{N}}^*$. We give in Section \ref{sect:properties} the properties of $(X_n,M_n)_{n\ge 0}$ and we determine necessary and sufficient conditions for existence and unicity of a stationary probability measure, in Section \ref{sect:invariant}. We would like to emphasize that $(X_n)$ is non-Markovian in general.

In Section~\ref{sec:contexttree}, we consider two particular cases of VLMC, associated with the simple infinite comb and the double infinite comb. In each of these two cases, the stationary measure can be explicitely calculated, in \cite{ccpp} for the simple comb and in the Appendix for the double comb. We make precise the correspondence between the process $(X_n,M_n)_{n\ge 0}$ defined in Section  \ref{sect:model} and the VLMC $(U_n)$ whose the last letter is $X_n$. Namely, we establish the dictionary between the stationary measure for the Markov chain $(X_n,M_n)_{n\ge 0}$ and the stationary measure for the VLMC (see Theorems \ref{prop:vlmc} and \ref{prop:vlmcd}). Thanks to these results, we do not have to worry about the point of view (couple letter/memory or VLMC)  when considering the persistent random walk
 $S_t:=\displaystyle\sum_{n=0}^t X_n$, under the stationary regime.
 
Section~\ref{sect:persistent} is devoted to the study of $(S_n)$. In particular, we determine  the explicit distribution of the r.v. $S_n$, see Proposition~\ref{prop:loicm}. Although the result is complicated, we are able to determine explicitely the generating function of the r.v. $S_{\tau+1}$, where $\tau$ is a geometric r.v. independent from $(X_n,M_n)_{n\ge 0}$. One way to compare the process $(S_n)$ with a classical random walk is to analyse how both processes fluctuate at infinity. We have the two following limit theorems, see Section~\ref{sec:largetime}:
\[
\lim_{n\to\infty}\frac{S_n}{n}=\Xi,\quad \mbox{and}\ \sqrt{n}\left( \frac{S_n}{n}-\Xi \right)\to \mathcal{N}(0,\sigma) 
\]
where $\mathcal{N}(0,\sigma)$ is a Gaussian distribution with $0$ mean and variance $\sigma^2$, $\Xi$ and $\sigma$ are constants which can be expressed in terms of the model parameters.

Finally in Section~\ref{sec:GITN} we study the persistent random walk. After a convenient scaling, its converges towards a ``generalized ITN'' as proved in Theorem~\ref{thm:convergenceGITN}. More precisely, we focus on the limit in law of Markov chains of the type $(X_n^\eps,M_n^\eps)$ which depends on a small parameter $\eps>0$. We suppose that $X^\eps_n$ takes its values in $\{-1,1\}$, $X_0^\eps=1$ and
\begin{align}
\label{eq:intro7}
\mathbb{P}\Big(X^\eps_{n+1}=1\Big| X^\eps_{n}=-1,\ M^\eps_n=k\Big)=f_1(k\eps)\eps+o(\eps)\\
\label{eq:intro8}
\mathbb{P}\Big(X^\eps_{n+1}=-1\Big| X^\eps_{n}=1,\ M^\eps_n=k\Big)=f_2(k\eps)\eps+o(\eps)
\end{align}
where $f_1,f_2:[0,\infty[\to\mathbb{R}$ are non negative and right continuous functions. Note that \eqref{eq:intro7} and \eqref{eq:intro8} mean that $(X^\eps_n)$ has a \emph{conservative} behaviour: if $X^\eps_n=1$ (resp. $X^\eps_n=-1$) the probability that $X^\eps_{n+1}$ changes, i.e. $X^\eps_{n+1}=-1$ (resp. $X^\eps_{n+1}=1$) is small for convenient $f_1, f_2$ and is measured by the parameter $\eps$.\\
Under additional assumptions, see the beginning of Section \ref{sec:GITN} for details, it is actually possible to rescale the triplet $(X^\eps_n,M^\eps_n,S^\eps_n)$ so that it converges as $\eps\to 0$. For simplicity, we only present the scaling procedure concerning $S_n^\eps$. The process $(S^\eps(t),\ t\ge 0)$ is piecewise linear and satisfies
\begin{align}
\label{eq:intro9}
S^\eps(n\eps)=\eps S_n^\eps, \quad \mbox{for any }\ n\in\mathbb{N}.
\end{align}
We prove (see Theorem \ref{thm:convergenceGITN} for a more complete result) that $(S^\eps(t),\ t\ge 0)$ converges in distribution, as $\eps\to 0$ to $(S^0(t),\ t\ge 0)$ where
\[
S^0(t)=\int_0^t (-1)^{N^0(s)}ds.
\]
Here, $(N^0(t))$ is the counting process with jump times $(\xi_n)_{n\ge 0}$:
\[
N^0(t)=\sum_{n\ge 0}\ind{\xi_n\le t},
\]
where $(\xi_{n+1}-\xi_{n},\ n\ge 0)$ is a sequence of independent r.v. such that $\xi_0=0$ and
\begin{align*}
\mathbb{P}(\xi_{2n+1}-\xi_{2n}\ge t)=\exp\left(-\int_0^t f_2(u)\, du\right)\\
\mathbb{P}(\xi_{2n+2}-\xi_{2n+1}\ge t)=\exp\left(-\int_0^t f_1(u)\, du\right)
\end{align*}
for any $t\ge 0$ and $n\ge 0$, where $f_1, f_2$ satisfies (\ref{eq:int-fini}).\\
The process $(S^0(t),\ t\ge 0)$ is called the \emph{Generalized Integrated Telegraph Noise} (see \cite{Herrmann-Vallois} for the ITN). It is both a semi-Markov process and a Piecewise Deterministic Markov Process \cite{Davis,Davisbook,Cocozza} and its trajectories look like a zig-zag.

%Notice that the results are presented for a general finite alphabet $\rond A := \{a_1,a_2,\ldots,a_K\} $ in Section~\ref{sect:model} and from Section~\ref{sec:contexttree}, $K=2$ and the letter process has two states $-1$ and $+1$. See Remark \ref{rem:gene:alphabet} for a generalization. Also notice that this study could be generalized as follows: consider, for any $a_i\in\rond A$, a (deterministic) sequence of real numbers $(\Gamma_{a_i}(k),\, k\ge 0)$. Set for $n\ge 0$,
%\begin{equation}
%\label{eq:def-persist-gen}
%S_{n+1}=S_n+\Gamma_{X_{n+1}}(n+1),
%\end{equation}
%with $S_0=\Gamma_{X_0}(0)$. We only deal in this paper with the case $\Gamma_{a_i}(t)=a_i$.

%%%%%%%%%%%%%%%%%%%%%%%%%%%%%%%%%%%%%%%%%%%%%%%%%%%%%%%%%%%%%%
%
%
%%%%%%%%%%%%%%%%%%%%%%%%%%%%%%%%%%%%%%%%%%%%%%%%%%%%%%%%%%%%%%
\mathversion{bold}
\section{The Markov chain $(X_n,M_n)_{n\geq 0}$}
\label{sect:model}
\setcounter{equation}{0}
\mathversion{normal}
\subsection{Definition}
Let us consider the finite set $\rond A = \{a_1,\ldots, a_{K}\}$ with $K>1$ elements. To each $a_i$ is associated a sequence $(\alpha_{i,n})_{n\ge 1}\in ]0,1[^{\mathbb{N}^ *}$ where $\mathbb{N}^*=\{1,2,3,\ldots\}$. We can now introduce the Markov chain $(X_n,M_n)_{n\in\mathbb{N}}$ valued in the state space $\{a_1,\ldots, a_{K}\}\times \mathbb{N}^*$ with transition probabilities
\begin{align}
\label{eq:transition}
\left\{\begin{array}{l}
Q\Big((a_i,n),(a_i,n+1)\Big)=1-\alpha_{i,n},\\
Q\Big((a_i,n),(a_j,1)\Big)=\alpha_{i,n}\ p_{i,j}, \quad 1 \le i\neq j\le K,\quad n\ge 1,
\end{array}\right.
\end{align}
where $\mathcal{P}:=(p_{i,j})$ is a given $K\times K$ transition matrix satisfying $p_{i,i}=0$ for all $1\le i\le K$, $p_{i,j}>0$ for all $i\neq j$ and $\sum_{j=1}^K p_{i,j}=1$ for all $i$. In fact, $p_{i,j}$ is the probability to move from $a_i$ to $a_j$ knowing that we leave $a_i$. 
\begin{figure}[h]
\centerline{\begin{tikzpicture}[scale=1.3]
\draw (0,0) -- (0,4) ;
\draw [>=stealth,->] (0,0) -- (0,1) ;
\draw (0,0) -- (6,0) ;
%\draw [>=stealth,->] (2,2) -- (2,1.5) ;
\draw (1,1) node[below]{{\scriptsize $(X_0,M_0)$}};
\draw (1,2) node[above]{{\scriptsize $(X_1,M_1)$}};
\draw (2,1) node[below]{{\scriptsize $(X_2,M_2)$}};
\draw (0,1) node[left]{$1$};
\draw (1,0) node[below]{$a_1$};
\draw (2,0) node[below]{$a_2$};
\draw (3,0) node[below]{$a_3$};
\draw (4,0) node[below]{$a_4$};
\draw (5,0) node[below]{$a_5$};
\draw [>=stealth,->](1,1) -- (1,1.5) ;
\draw [>=stealth,->](1,2) -- (1.5,1.5) ;
\draw (1,1) -- (1,2) node[midway, above, sloped]{{\scriptsize $1-\alpha_{1,1}$}};
\draw (1,2) -- (2,1) node[midway, above, sloped]{{\scriptsize $\alpha_{1,2}$}} ;
\draw [>=stealth,->] (2,1) -- (2,1.5) ;
\draw [>=stealth,->] (2,2) -- (2,2.5) ;
\draw [>=stealth,->] (2,3) -- (2,3.5) ;
\draw [>=stealth,->] (2,4) -- (2.5,3.25) ;
\draw [>=stealth,->] (4,1) -- (4,1.5) ;
\draw [>=stealth,->] (4,2) -- (3.5,1.5) ;
\draw [>=stealth,->] (3,1) -- (3,1.5) ;
\draw (2,1) -- (2,2) node[midway, below, sloped]{{\scriptsize $1-\alpha_{2,1}$}};;
\draw (2,2) -- (2,3) node[midway, above, sloped]{{\scriptsize $1-\alpha_{2,2}$}};;
\draw (2,3) -- (2,4) node[midway, above, sloped]{{\scriptsize $1-\alpha_{2,3}$}};;
\draw (2,4) -- (4,1) node[midway, above, sloped]{{\scriptsize $\qquad\alpha_{2,4}$}} ;
\draw (4,1) -- (4,2) ;
\draw (4,2) -- (3,1) ;
\draw (3,1) -- (3,3) ;
\draw [dashed, thick] (3,3) -- (3,4) ;
\draw [dotted, very thin, gray] (0,0) grid (6,4);
\draw (0,0) node[left]{$0$} ;
\draw (6,0) node[right]{$X_n$} ;
\draw (0,4) node[left]{$M_n$} ;
\draw (1,1) node {$\bullet$};
\draw (1,2) node {$\bullet$};
\draw (2,1) node {$\bullet$};
\draw (2,2) node {$\bullet$};
\draw (2,3) node {$\bullet$};
\draw (2,4) node {$\bullet$};
\draw (3,1) node {$\bullet$};
\draw (3,2) node {$\bullet$};
\draw (3,3) node {$\bullet$};
\draw (4,1) node {$\bullet$};
\draw (4,2) node {$\bullet$};
\draw (3,4) node {$\bullet$};
\end{tikzpicture}}
\caption{\label{fig-path}A path description of the process $(X_n,M_n)_{n\ge 0}$}
\end{figure}

\noindent Moreover, in order to deal with VLMC later on, we extend the definition of the Markov chain to the state space $\{a_1,\ldots, a_{K}\}\times\overline{\mathbb{N}}^*$ with $\overline{\mathbb{N}}^*=\mathbb{N}^*\cup\{\infty\}$. Therefore we introduce $\alpha_{i,\infty}\in]0,1[$ for all $1\le i\le K$ such that
\begin{equation}
 \label{eq:extension}
\left\{\begin{array}{l}
Q\Big((a_i,\infty),(a_i,\infty) \Big)=1-\alpha_{i,\infty}\\
Q\Big((a_i,\infty),(a_j,1) \Big)=\alpha_{i,\infty}\ p_{i,j},\quad i\neq j.
\end{array}\right.
\end{equation}
Note that $\alpha_{i,k}$ is the probability of changing letter after a run of length $k$ of $a_i$, that is
\begin{equation}\label{desc-alpha}
\alpha_{i,k}=\mathbb{P}(X_{n+1}\neq a_i\vert X_n=a_i,\ M_n=k).
\end{equation}
There are strong links between $(X_n)$ and $(M_n)$. In particular,
if $M_0=1$, $M_n$ can be expressed with $X_0,\ldots,X_n$. Indeed, if the sequence $(X_j)_{j=0, \ldots ,n}$ is constant then $M_n=n+1$ and $M_n=\inf\{1\le i\le n;\ X_{n-i}\neq X_n\}$ otherwise. In other words, one has
\begin{align}
\label{eq:def:M}
M_n&=1+\sup\{0\le i\le n,\  X_{n-j}=X_n, \ \forall j\in\{ 0,\ldots,i \}\}\\
&=\inf\{ 0\le i \le n,\ X_{n-i}\neq X_n \}.\nonumber
\end{align}
Let us explain how moves $(X_n,M_n)$ in the case $M_0=1$ and $X_0=a_i$. The variable $M_n$ increases by one
unit at each time until $X_n$ switches to $a_j\neq a_i$. At that
first jump time, the memory is reset to $1$ and so on... So that $M_n$
represents the variable memory of $(X_t)_{0\le t\le n}$ since it
counts the last consecutive stays (at $X_n$) before $n$. Moreover
the dynamics of the jumps of $X_{n}$ is governed by the value of
$M_n$. In Figure \ref{fig-path}, we have drawn the following trajectory of $(X_n,M_n)$ corresponding to the values $a_1,a_1,a_2,a_2,a_2,a_2,a_4,a_4,a_3,a_3,a_3,a_3 \ldots $ of $X_n$.

\noindent Let us note that in the particular case: $K=2$,  $a_1=0$, $a_2=1$ and $\alpha_{j,n}=\alpha_j$
for all $n\ge 1$, then
$(X_n)_{n\ge 0}$ is a sequence of independent Bernoulli random
variables. \mathversion{bold}
\subsection{Properties of the Markov chain $(X,M)$}
\label{sect:properties}

\mathversion{normal} First we investigate under which conditions either $(X_n)$ or $(M_n)$ is Markov. Secondly we prove existence of invariant probability measure and
finally we present a path description for the process
$(X_n)$.
\subsubsection{Link between the margins}
A natural question arising about a 2-dimensional Markov chain $(X_n,M_n)_{n\in \Nset}$ is to know whether the margins are Markovian too. The following proposition says that in general case, neither $(X_n)_{n\in\mathbb{N}}$ nor $(M_n)_{n\in\mathbb{N}}$ is a Markov chain.
\begin{prop}
\label{prop:markov-prop} Assume that $M_0=1$.
\begin{enumerate}[(i)]
\item The margin process $(X_n)_{n\in\mathbb{N}}$ is Markovian if and only if for all $1\le i\le K$, $n\mapsto \alpha_{i,n}$ is constant. In that case the transition matrix $Q^X$ of $X$ is given by:
\begin{align}\label{pi-dedud}
Q^X(i,j)=\left\{\begin{array}{ll} 1-\alpha_{i,1} & \mbox{if}\
j=i\\
 \alpha_{i,1}p_{i,j}\ & \mbox{if}\ j\neq i.
\end{array}
\right.
\end{align}
\item The margin process $(M_n)_{n\in\mathbb{N}}$ is Markovian for any initial condition $X_0$ if and only if for all $n\geq 1$, the function $i\mapsto \alpha_{i,n}$ is constant. In that case, 
%\begin{equation}
%\label{eq:equal}
%\alpha_{i,k}=\alpha_{1,k},\quad\forall k\ge 1, \forall i\in\{1,\ldots,K\}.
%\end{equation}
%Under \eqref{eq:equal}, 
the transition matrix $Q^M$ of $M$ is
\[
Q^M(n,j)=\left\{\begin{array}{ll} 1-\alpha_{1,n} & \mbox{if}\
j=n+1\\
\alpha_{1,n} &\mbox{if}\ j=1.
\end{array} \right.
\]
%3) Let $X_0=a_i$. Then the sigma-algebra
%$\sigma\{X_0,X_1,\ldots,X_n\}$ generated by $X_0$,..., $X_n$
%contains $\sigma\{M_0,M_1,\ldots,M_n\}$. Moreover in the
%particular case $K=2$, these sigma-algebras are equal together.
\end{enumerate}
\end{prop}
\pff 
\begin{enumerate}[(i)]
\item For a given vector $(x_0, \ldots, x_n)\in \{a_1, \ldots, a_K\}^{n+1}$, let us first denote
\[
\delta_{i,n}:=\mathbb{P}(X_{n+1}=a_i\vert X_n=x_n,\ldots,X_0=x_0).
\]
According to \eqref{eq:def:M} let us introduce:
\begin{equation}\label{eq:def:m}
m_n=1+\sup\{0\le i\le n:\ x_{n-j}=x_n,  \ \forall j\in\{ 0,\ldots,i \}\}.
\end{equation}
We have to distinguish two cases.
\begin{enumerate}
\item If $x_n=a_i$ then $M_n=m_n$ and therefore $\delta_{i,n}=1-\alpha_{i,m_n}$. We can choose different values of $x_2$,...,$x_{n-1}$ such that $m_n=1,2,\ldots,n$. Hence if $(X_{n})$ is Markovian then $\delta_{i,n}$ is independent of $n$ and $(x_0, \ldots, x_{n-1})$ and $(\alpha_{i,k})_{k\ge 1}$ is constant.
\item If $x_n=a_j\neq a_i$ then $\delta_{i,n}=\alpha_{j,m_n}\ p_{j,i}=\alpha_{j,1}\ p_{j,i}$ implying that $(X_n)$ is effectively Markovian.
\end{enumerate}
\item Let us study the process $(M_n)_{n\ge 0}$. Set
\[
d_{i,n}:=\mathbb{P}(X_0=i,\, M_0=1,\, M_1=2,\ldots, M_n=n+1,\,
M_{n+1}=1).
\]
We have
\[
d_{i,n}=\mathbb{P}(X_0=i,\, M_0=1,\, X_1=i,\,
M_1=2,\ldots,X_n=i,\, M_n=n+1,\, M_{n+1}=1).
\]
Since $(X_n,M_n)$ is a Markov chain, using \eqref{eq:transition}
and \eqref{desc-alpha} we get
\[
d_{i,n}=(1-\alpha_{i,1})\times\ldots\times(1-\alpha_{i,n})\alpha_{i,n+1}.
\]
Suppose that $(M_n)$ is a Markov chain, with transition matrix
$Q^M$. Then
\begin{align*}
d_{i,n}&=Q^M(n+1,1)\mathbb{P}(X_0=i,\,M_0=1,\,\ldots,\,
M_n=n+1)\\
&= Q^M(n+1,1)(1-\alpha_{i,1})\times\ldots \times(1-\alpha_{i,n}).
\end{align*}
Consequently, $\alpha_{i,n+1}=Q^M(n+1,1)$ in independent of $i$ and thus $\alpha_{i,n+1}=\alpha_{1,n+1}$ for all $i$ and $n$.\\
As for the converse, since $i\mapsto \alpha_{i,n}$ is constant, it is clear that
\eqref{eq:transition} implies that $(M_n)$ is Markov. \QED
\end{enumerate}
%\begin{rem}\label{rem:filtra}\penote{enlever cette remarque ?} Suppose that $M_0=1$ and $X_0=a_i$.
%It is clear that \eqref{eq:def:M} implies that
%$\sigma(M_0,\ldots,M_n)\subset \sigma(X_0,\ldots,X_n)$. In the
%particular case $K=2$, let us introduce
%$R_n=\sum_{k=1}^n\ind{M_k=1}$. $R_n$ is the number of jumps of
%$X$ up to time $n$. Consider $r_n$ the rest of $R_n$ by $2$:
%$R_n=2q_n+r_n$ where $r_n=0$ or $1$. Then $X_n=a_{r_n+1}$. In
%particular $\sigma(M_0,\ldots,M_n)=\sigma(X_0,\ldots,X_n)$.
%\end{rem}
\begin{rem} The one-dimensional memory process $(M_n)_{n \in \Nset}$ could be replaced by a $K$-dimensional process. For each state $a_i$, define
\[
\mathcal{M}_n^{(a_i)}=\inf\{0\le k\le n;\ X_{n-k}\neq a_i\}.
\]
There is a one-to-one correspondence between $(X_n,M_n)$ and
$(\mathcal{M}_n^{(a_1)},\ldots,\mathcal{M}_n^{(a_K)})$. Consequently the vector memory $(\mathcal{M}_n^{(a_1)},\ldots,\mathcal{M}_n^{(a_K)})$ is a Markov chain. For instance $(X_n,M_n)=(a_3,4)$
corresponds to
$(\mathcal{M}_n^{(a_1)},\ldots,\mathcal{M}_n^{(a_K)})=(0,0,4,0,\ldots,0)$.
Indeed if the $k$\textsuperscript{th} coordinate of the vector does not vanish then
$X_n=a_k$. This permits to recover $X_n$ via $(\mathcal{M}_n^{(a_1)},\ldots,\mathcal{M}_n^{(a_K)})$, as for $M_n$, we have $M_n=\mathcal{M}_n^{X_n}$.\end{rem}
%Let us define $M_n^{(i)}$ the memory corresponding to state $a_i$:
%\[
%M_n^{(i)}=\inf\{0\le k\le n;\ X_{n-k}\neq a__Ã§i\}.
%\]
%As shows Proposition \ref{prop:markov-prop}, $(M_n)$ is not in general Markovian, however $(M^{(1)}_n,\ldots,M^{(K)}_n)$ is a Markov chain.\\[5pt]
\mathversion{bold}
\subsubsection{Invariant probability measure for $(X_n,M_n)_{n\geq 0}$} 
\label{sect:invariant}
\mathversion{normal}
Let us now
investigate the existence of an invariant probability measure. It
is convenient to introduce for all $1\le i\le K$,
\begin{equation}
 \label{eq:cond-gene}
\Theta_i:=\sum_{n\ge 1}\prod_{k=1}^{n-1} (1-\alpha_{i,k}),
\end{equation}
and for $m\geq 1$,
\begin{equation}
\label{eq:}
\mathcal{P}_i(m):=\prod_{k=1}^{m-1}(1-\alpha_{i,k}),
\end{equation}
with the convention $\prod_1^0=1$.\\
$\mathcal{P}_i(k+1)$ represents the conditional probability that the 
process $(X_n)$ stays at least a time interval of length $k$ in the same state $i$ 
\[\mathcal{P}_i(k+1)=\PP\left(X_1= \ldots= X_{k}=i\big|X_0=i,M_0=1\right).\]

%Recall that if $\mathcal{P}$ is an irreducible transition matrix,
%Frobenius' theorem (see, for instance Seneta \cite{Seneta}) \penote{a deplacer} ensures the
%existence of a unique invariant probability
%$v^*=(v_1^*,\ldots,v^*_K)$ satisfying
%\begin{equation}
%\label{eq:frobe}
%v^*=v^\intercal \mathcal{P}^*
%\end{equation}
%\bnote{ super mal dit. qui est v, P* ?} where $b^\intercal$ stands for the transpose of the matrix $b$.\\
%In other words, $(v^*)^\intercal$ is the unique invariant probability for the transition matrix $\mathcal{P}$.
\begin{prop}
 \label{prop:stat-gen}
 
 Let us denote $\mathcal{P}=(p_{i,j})$ a given irreducible transition matrix.
 % defined by the transition probabilities (\ref{eq:transition}). 
 \begin{enumerate}[(i)]
 \item Then the Markov chain $(X_n,M_n)_{n\ge 0}$ with transition probabilities defined by
 \eqref{eq:transition} and \eqref{eq:extension} admits a invariant probability
 measure $\nu$ on the space $\{a_1,\ldots, a_{K}\}\times \overline{\mathbb{N}}^*$ 
 if and only if $\Theta_1$,...,$\Theta_K$ defined by \eqref{eq:cond-gene} are all finite. This invariant probability measure is unique.
\item Moreover, if we denote by $v^*$ the unique positive vector associated with the largest eigenvalue of $\mathcal{P}=(p_{i,j})$ by Frobenius's theorem, then
$\nu(a_i,\infty)=0$ for all $1\le i\le K$ and $n\geq 1$,
\[\nu(a_i,n)=\frac{v_i^*}{\langle \Theta, v^*\rangle} \mathcal{P}_i(n)  \]
where $\Theta={}^t(\Theta_1,\ldots,\Theta_K)$ and $\langle \Theta,v^*\rangle=\sum_{i=1}^K \Theta_i v^*_i$.
\end{enumerate}
\end{prop}
\begin{rem}
The invariant measure $\nu$ can be decomposed in the following way: for $1\le i\le K$ and $n\ge 1$,
\begin{equation}\label{eq:decom-sup}
 \nu(a_i,n)=\nu^{X}(a_i)\nu_i(n), 
 \end{equation}
where
\[
 \nu^X(a_i)=\frac{\Theta_i v^*_i}{\langle \Theta,v^*\rangle}\quad\mbox{and}\quad \nu_i(n)=\frac{\mathcal{P}_i(n)}{\Theta_i}.
\]
If $(X_0,M_0)\sim \nu$, then, for any $n\ge 1$, $\nu^X$ is the law of $X_n$, and $\nu_i$ is the conditional distribution of $M_n$, given $X_n=i$.\\
 Let us consider the particular case when for all $1\le i\le K$ and $n\ge 1$,
\[
 1-\alpha_{i,n}=\frac{\rho_i}{n},\quad \mbox{with }\rho_i>0.
\]
After straightforward calculations, we obtain $\Theta_i=e^{\rho_i}$ and
\[
 \nu^X(a_i)=\frac{v_i^* e^{\rho_i}}{\langle \Theta, v^*\rangle},\qquad \nu_i(n)=\frac{\rho_i^{n-1}}{(n-1)!}\, e^{-\rho_i}.
\]
In other words, if $(X_0,M_0)\sim \nu$ then the distribution of the couple $(X_n,M_n)$ can be described as follows: $X_n$ is chosen first with the probability $\nu^X$ and afterwards, conditionally on $X_n=a_i$, $M_n$ is Poisson distributed with parameter $\rho_i$.
\end{rem}
\noindent {\sc Proof of Proposition \ref{prop:stat-gen}.} 

For notational simplicity, we shall fix $a_i=i$ for all $1\le i\le K$.
\paragraph{Step 1 --- Invariant measure:}

Let $\nu$ be a non-negative measure. Since $(X_n,M_n)$ is valued in the state space $\{1,2,3,\ldots,K\}\times\overline{\mathbb{N}}^*$, $\nu$ is an invariant measure if and only if
\begin{align}
\label{equat1}
\nu(i,k)&=\sum_{\ell\ge1}\Big\{ \nu(i,\ell)Q\Big( (i,\ell),(i,k) \Big)+\sum_{j\neq i} \nu(j,\ell)Q\Big( (j,\ell),(i,k) \Big)\Big\}\nonumber\\
&=\nu(i,k-1)(1-\alpha_{i,k-1})\ind{k>1}+\ind{k=1}\sum_{j\neq i}p_{j,i}\sum_{\ell\ge 1}\nu(j,\ell)\alpha_{j,\ell},
\end{align}
for any $1\le i \le K$, and
\begin{equation}\label{eq:nullite}
 \nu(i,\infty)=\nu(i,\infty)(1-\alpha_{i,\infty}).
\end{equation}
Obviously \eqref{eq:nullite} implies that $\nu(i,\infty)=0$ for all $1\le i\le K$.
Relation \eqref{equat1}, with $k\ge 2$ is equivalent to $\nu(i,k)=\nu(i,k-1)(1-\alpha_{i,k-1})$ which implies for $k\geq 2$
\begin{equation}
\label{eq:interm1}
\nu(i,k)=\nu(i,1)\prod_{r=1}^{k-1}(1-\alpha_{i,r}) = \nu(i,1) \rond P_i(k).
\end{equation}
The particular situation $k=1$ in \eqref{equat1} and \eqref{eq:interm1} leads to
\begin{align}\label{eq:ajout}
\nu(i,1)%&=\sum_{j\neq i}p_{j,i}\sum_{l\ge 1}\alpha_{j,l}\nu(j,l)\nonumber\\
&=\sum_{j\neq i}p_{j,i}\Big(\sum_{\ell\ge 1}\alpha_{j,\ell}\rond P_j(\ell)\Big)\nu(j,1).
\end{align}
Using \eqref{eq:interm1} and \eqref{eq:cond-gene} we get:
\[
\sum_{1\le i\le K,\ n\ge 1}\nu(i,n)=\sum_{i=1}^K\nu(i,1)\sum_{n\ge
1}\rond P_i(n)=\sum_{i\ge 1}\nu(i,1)\Theta_i.
\]
Finally $\nu$ is a probability measure iff $\nu(i,k)$ is given by \eqref{eq:interm1} for any $k\ge 1$,
the vector ${}^t(\nu(1,1),\ldots,\nu(K,1))$ solves \eqref{eq:ajout} and
\begin{equation}
\label{eq:somme}
\sum_{i=1}^K \nu(i,1)\Theta_i=1.
\end{equation}
\paragraph{Step 2 --- Necessary condition:}

Assume that
\begin{equation}\label{eq:repet}
 \Theta_i<\infty,\quad\forall i\in\{1,\ldots,K\}.
\end{equation}
Writing $\alpha_{j,l}=-(1-\alpha_{j,l})+1$, we develop the expression \eqref{eq:ajout} using \eqref{eq:repet}:
\begin{align*}
 \nu(i,1)&=\sum_{j\neq i}p_{j,i}\Big\{ \sum_{\ell\ge 1} \rond P_j(\ell)-\rond P_j(\ell+1)\Big\}\ \nu(j,1)\\
&=\sum_{j\neq i}p_{j,i}\ \nu(j,1).
\end{align*}
The vector $v:={}^t(\nu(1,1),\nu(2,1),\ldots,\nu(K,1))$ satisfies
\begin{equation}
 \label{eq:vect-transp}
v={}^t\mathcal{P} v,\quad \mbox{with}\ \mathcal{P}=(p_{i,j}).
\end{equation}
Let $v^*$ be the unique positive vector associated with the largest eigenvalue of $\mathcal{P}=(p_{i,j})$ by Frobenius' theorem, then there exists $\lambda>0$ such that
%As a result, ${}^t(\nu(1,1),\ldots,\nu(K,1))$ solves \eqref{eq:ajout} if and only if there exists $\lambda>0$ such that
\[
 (\nu(1,1),\nu(2,1),\ldots,\nu(K,1))=\lambda {}^t v^*.
\]
Using \eqref{eq:somme} we deduce:
\begin{align*}
\sum_{i=1}^K\nu(i,1)\Theta_i=\lambda \sum_{i=1}^K\Theta_i v^*_i=\lambda\langle \Theta,v^*\rangle.
\end{align*}
Hence $\lambda=1/\langle \Theta,v^*\rangle$ and by \eqref{eq:interm1}, $\nu$ is determined by $\displaystyle\nu(i,n)=\frac{v_i^*}{\langle \Theta, v^*\rangle} \mathcal{P}_i(n)$, which gives existence and unicity of $\nu$.\\
\paragraph{Step 3 --- Sufficient condition:}
Conversely let us assume the existence of an invariant
probability measure $\nu$. We shall prove \eqref{eq:repet}.
Obviously \eqref{eq:ajout} implies that if $\nu(i,1)=0$ for some
$i$, then $\nu(j,1)=0$ for all $j$. Therefore $\nu=0$ which
contradicts the fact that $\nu$ is a probability measure. Hence
$\nu(i,1)>0$ for all $1\le i\le K$. It is clear that
\eqref{eq:somme} implies \eqref{eq:repet}.\QED

\noindent Since the Markov chain $(X_n,M_n)_{n\in\mathbb{N}}$ admits an invariant probability measure, we can extend its definition to $\mathbb{Z}$ (instead of $\mathbb{N}$) such that it is stationary. This extension will be usefull to connect with certain Variable Length Markov Chains
(defined later in Section~\ref{sec:contexttree}).
% it is convenient to introduce a construction of the
%above Markov chain on the canonical space where the time space is $\mathbb{Z}$
%instead of $\mathbb{N}$. Let us recall the classical procedure.
%Let $(Z_n)_{n\ge 0}$ be a $E$-valued Markov chain ($E$ is assumed
%to be discrete) which admits a stationary probability measure
%$\nu$. We define on the canonical space $\Omega_c=E^\mathbb{Z}$, a
%family of probability measures depending on the parameters
%$i\in\mathbb{Z}$ and $n\ge 0$:
%\[
% \mu^{i,n}(A_i\times\ldots\times A_{i+n})=\mathbb{P}_\nu(Z_0\in A_i,\ldots, Z_n\in A_{i+n}),
%\]
%where $\mathbb{P}_\nu$ is the probability under which $(Z_n)_{n\ge 0}$ is stationary.
%Since $(Z_n)_{n\ge 0}$ is $\mathbb{P}_\nu$-invariant, the family
%$\{\mu^{i,n};\, i\in\mathbb{Z},\, n\ge 0\}$ is consistent. We
%deduce the existence of a probability measure on $\Omega_c$
%denoted $\mathbb{Q}$ satisfying
%\[
%\mathbb{Q}(\xi_i\in A_i,\ldots, \xi_{i+n}\in A_{i+n})
%=\mathbb{P}_\nu(Z_0\in A_i,\ldots, Z_n\in A_{i+n}),
%\]
%where $(\xi_n)_{n\in\mathbb{Z}}$ is the canonical process. By this procedure, we have just construct a process $(\xi_n)_{n\in\mathbb{Z}}$ associated with the dynamics of $(Z_n)_{n\ge 0}$ under $\mathbb{P}_\nu$.\\[5pt]
%These arguments can of course be applied to the Markov chain
%$(X_n,M_n)_{n\ge 0}$. For notational simplicity, the canonical
%process shall be also denoted
%$(X_n,M_n)_{n\in\mathbb{Z}}$ in the following.
\begin{rem}
 \label{rem:ajout}
Since $\nu$ is the invariant probability measure of the Markov chain $(X_n,M_n)_{n\in\mathbb{Z}}$ then $(X_{-n},M_{-n})_{n\in\mathbb{Z}}$ is a Markov chain with invariant probability measure $\nu$ and transition probabilities $\widehat{Q}$ where:
\[
 \nu(x)Q(x,y)=\nu(y)\widehat{Q}(y,x),\quad \forall x,y\in\{a_1,\ldots,a_K\}\times\mathbb{N}^*.
\]
From \eqref{eq:transition} and Proposition~\ref{prop:stat-gen} we easily obtain
\begin{align*}
 \left\{\begin{array}{l}
\widehat{Q}\left( (a_i,n+1),(a_i,n) \right)=1,\quad \forall i\in\{ 1,\ldots,K \},\ n\ge  1,\\[8pt]
\displaystyle \widehat{Q}\left( (a_j,1),(a_i,n) \right)=\frac{v^*_i}{v^*_j}\, p_{i,j}\alpha_{i,n}\mathcal{P}_i(n),\quad i\neq j, \ n\ge 1.
\end{array}
 \right.
\end{align*}
\end{rem}
%%%%%%%%%%%%%%%%%%%%%%%%%%%%%%%%%%%%%%%%%%%%%%%%%%%%%%%%%%%%%%%%%%%%%%%%%%%%%%%%%%
\mathversion{bold}
\subsubsection{Paths description of $X$}
\mathversion{normal}
From now on, for notational simplicity, we only consider the case $K=2$. The trajectory $n\mapsto X_n$ is determined as soon as the transition times between the different
states are known. Let us define $T_0=0$ and the sequence of
stopping times for $n\ge 1$,
\begin{equation}
\label{eq:changt} T_n=\inf\Big\{ i\ge T_{n-1}:\ X_i\neq
X_{T_{n-1}} \Big\}.
\end{equation}
\begin{prop}
\label{prop:changtpente-sym} 
\begin{enumerate}[(i)]
\item Let us assume that $\Theta_1$ and $\Theta_2$ defined by \eqref{eq:cond-gene} are finite. Then the random variables
$(T_{n+1}-T_n)_{n\ge 1}$ are almost surely finite and independent.
\item 
\begin{enumerate}
\item If $X_0=a_2$ and $M_0=m\ge 1$. Then for all $i\ge 1$ and $n\ge 1$,
\begin{equation}
\label{eq:dist+}
\mathbb{P}(T_{1}=i)=\alpha_{2,m+i-1}\prod_{j=m}^{m+i-2}(1-\alpha_{2,j}),
\end{equation}
and
\begin{align}
\label{eq:dist++}
&\mathbb{P}(T_{2n+1}-T_{2n}=i)=\alpha_{2,i}\mathcal{P}_2(i),\\
\label{eq:dist+-+}
&\mathbb{P}(T_{2n}-T_{2n-1}=i)=\alpha_{1,i}\mathcal{P}_1(i).
\end{align}
\item If $X_0=a_1$ and $M_0=m\geq 1$ then \eqref{eq:dist+} and \eqref{eq:dist++} (resp.
\eqref{eq:dist+-+}) are still valid after replacing
$(\alpha_{2,\bullet})$ by $(\alpha_{1,\bullet})$ (resp.
$(\alpha_{1,\bullet})$ by $(\alpha_{2,\bullet})$).
\end{enumerate}
\end{enumerate}
\end{prop}
\begin{rem}\label{rem:2.6}
\begin{enumerate}
\item Note that, if $X_0=a_2$ and $M_0=m$ then for all $n\geq 1$,
\[
\mathbb{P}(T_1\ge
i)=\prod_{j=m}^{m+i-2}(1-\alpha_{2,j}),\quad\mathbb{P}(T_{2n+1}-T_{2n}\ge
i)
%=\prod_{j=1}^{i-1}(1-\alpha_{2,j})
 = \mathcal{P}_2(i)
\]
and
\(
\mathbb{P}(T_{2n}-T_{2n-1}\ge
i)
%=\prod_{j=1}^{i-1}(1-\alpha_{1,j}) 
= \mathcal{P}_1(i).
\)
%\item Suppose\bnote{ca sert plus tard ? A enlever ?} that $X_0=a_2$ and there exists $k\ge 1$ such that $\alpha_{2,k}=1$ (resp. $\alpha_{1,k}=1$) then $T_{2n+1}-T_{2n}\le k$ (resp. $T_{2n}-T_{2n-1}\le k$) for $n\ge 1$. Similarly if $\alpha_{2,k}=0$ (resp. $\alpha_{1,k}=0$) for some $k\ge 1$, then $T_{2n+1}-T_{2n}\neq k$ (resp. $T_{2n}-T_{2n-1}\neq k$) for any $n\ge 1$. \\
\item Between two consecutive jump times, the memory increases linearly
\begin{align}
\label{eq:R1}
M_{T_n+t}=1+t,\quad 0\le t< T_{n+1}-T_n,\ n\ge 1.
\end{align}
\end{enumerate}
\end{rem}
\noindent {\sc{Proof of Proposition \ref{prop:changtpente-sym}}} 

\noindent Let us consider $X_0=a_2$ and $M_0=m$. Then
\[
\mathbb{P}(T_1=i)=\mathbb{P}(X_1=a_2,X_2=a_2,\ldots,X_{i-1}=a_2,X_i=a_1).
\]
Using the Markov property, we deduce
\begin{align*}
\mathbb{P}(T_1=i)=\prod_{j=m}^{m+i-2}Q\left(
(a_2,j),(a_2,j+1) \right)Q\left( (a_2,m+i-1),(a_1,1)\right).
\end{align*}
Equation~\eqref{eq:dist+} is therefore a direct consequence of
\eqref{eq:transition}.\\
Using $\alpha_{2,m+i-1}=1-(1- \alpha_{2,m+i-1})$, it is easy to
deduce that
\[
\sum_{i\ge
1}\alpha_{2,m+i-1}\prod_{j=m}^{m+i-2}(1-\alpha_{2,j})=1.
\]
This shows that $\mathbb{P}(T_1<\infty)=1$.

\noindent Morever, conditioning by $X_0=a_2$ and $M_0=m$, for $j\geq 1$ and $i\geq 1$, one has
\begin{align*}
\PP\left(T_2-T_1=j, T_1=i\right)&=\PP\left(X_1=\ldots=X_{i-1}=a_2, X_i=\ldots=X_{j+i-1}=a_1, X_{j+i}=a_2\right)\\
&= \alpha_{2,m+i-1}\prod_{j=m}^{m+i-2}(1-\alpha_{2,j})\prod_{\ell=1}^{j-1}(1-\alpha_{1,\ell})\alpha_{1,j}\\
&= \PP\left(T_1=i\right)\PP\left(T_2-T_1=j\right),
\end{align*}
which leads to the independence between $T_1$ and $T_2-T_1$. The independence of $T_3-T_2$ and $(T_1, T_2-T_1)$ can be proved similarly. The proof of (ii) (a) of Proposition \ref{prop:changtpente-sym} follows by induction. The proof for (ii) (b) is analog. \QED
%Since $T_1$ is a finite stopping time, $X_{T_1}=a_1$, $M_{T_1}=1$,
%then under $\{X_0=a_2,\ M_0=k\}$, $(X_{t+T_1},M_{t+T_1})_{t\ge 0}$
%is distributed as $(X_t,m(t))$ under $\{X_0=a_1,\ M_0=1\}$. This
%proves \eqref{eq:dist+-+} with $n=1$. \QED

% In order to compute
%$T_{2n+1}-T_{2n}$ for $n\ge 0$, it suffices to note that, due to
%the definition of the stopping times, $X_{T_{2n}}=1$ and
%$M_{T_{2n}}=1$. Let us define
%\[
%\psi(i,l,2n+1):=\mathbb{P}_{(1,k)}(T_{2n+1}-T_{2n}=i\vert
%T_{2n}=l).
%\]
%Then
%\begin{align*}
%\psi(i,l,2n+1)&=\mathbb{P}_{(1,k)}(Y_{l+1}=1,\ldots, Y_{l+i-1}=1,Y_{l+i}=-1\vert T_{2n}=l)\\
%&=\prod_{j=1}^{i-1}\pi\Big( (1,j),(1,j+1) \Big)\pi\Big( (1,i),(-1,1) \Big)\\
%&=\alpha_{2,i}\prod_{j=1}^{i-1}(1-\alpha_{2,j}).
%\end{align*}
%Let us note that this expression does not depend on the variables $l$ and $n$.\\
%At time $T_{2n-1}$ the $2$-dimensional Markov chain $(X,M)$ is
%equal to $(-1,1)$. Hence the arguments used to prove
%\eqref{eq:dist+-+} are similar to those used just before.

%%%%%%%%%%%%%%%%%%%%%%%%%%%%%%%%%%%%%%%%%%%%%%%%%%%%%%%%%%%%%%%%%%
%
%
%%%%%%%%%%%%%%%%%%%%%%%%%%%%%%%%%%%%%%%%%%%%%%%%%%%%%%%%%%%%%%%%%%
\mathversion{bold}
\section{The variable length Markov Chain $(U_n)_{n \geq 0}$}
\mathversion{normal}
\label{sec:contexttree} 
\setcounter{equation}{0}

In this section, the relation between the Markov chain $(X_n, M_n)_{n \in \mathbb{Z}}$ valued in $\{0,1\}\times\overline{\mathbb{N}}^*$ and the VLMC $(U_n)_{n\geq 0}$ introduced in Section~\ref{intro} is highlighted by the Theorems~\ref{prop:vlmc} and \ref{prop:vlmcd}. For two very particular variable length Markov chains, we prove that these two models are equivalent.
%We begin with the definition of a general variable length Markov chain. The following presentation comes from \cite{ccpp}. Let $\mathcal{L}= \{0,1\}^{-\mathbb{N}}$ be the set of left-infinite words. Consider a complete (each node has $0$ or $2$ children) binary tree whose finite leaves $\mathcal{C}$ are words. Each leaf $c$ (not necessarily finite) is labelled with a Bernoulli distribution, respectively denoted by $q_{c}$. This probabilized tree is called a \emph{context tree}. The VLMC (Variable Length Markov Chain) associated with a context tree is the $\mathcal{L}$-valued Markov chain $(U_n)_{n\geqslant 0}$ defined by the transitions
%\begin{equation}
% \label{eq:def:VLMC}
 %\PP(U_{n+1} = U_n\alpha | U_n)= q_{\petitlpref (U_n)}(\alpha)
% \end{equation}
%where $\alpha\in \{0,1\}$ is any letter and $\lpref(U_n)$ denotes the first suffix of $U_n$ (reading from right to left) appearing as a leaf of the context tree. Notice that the VLMC is entirely determined by the data $q_{c}, c \in \mathcal{C}$. Moreover the order of dependence from the past depends on the past itself.
We consider two cases of VLMC for two specific context trees: the simple infinite comb and the double infinite comb. %We compare these models with the model defined in Section~\ref{sect:model}.

\noindent From now on and until the end of this paper, for the sake of simplicity, we only consider the case $K=2$.

\subsection{The simple infinite comb}
\label{simplecomb}
Let us consider the alphabet $\{a_1,a_2\}$ with $a_1=0$ and $a_2=1$. We associate with a Markov chain of type $(X_n,M_n)_{n\in\g Z}$ defined in Section \ref{sect:model} a unique VLMC and vice versa. Abusing words, this VLMC is called the \emph{infinite comb}. We refer to \cite{ccpp} for a complete definition.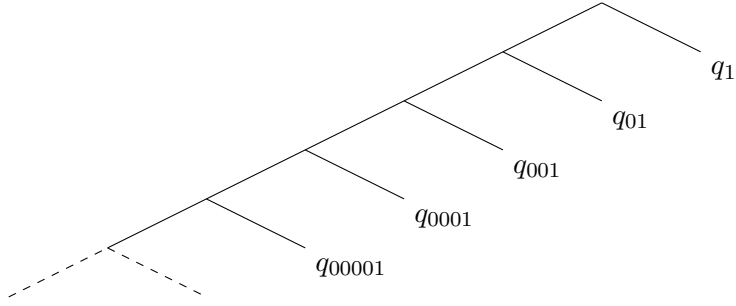
\begin{figure}
\centerline{\begin{tikzpicture}[scale=1.3]
\draw (1,0.5) -- (6,3) ;
%\draw [>=stealth,->] (0,0) -- (0,1) ;
\draw [dashed] (0,0) -- (1,0.5) ;
\draw [dashed] (1,0.5) -- (2,0);
\draw (6,3) -- (7,2.5);
\draw (5,2.5) -- (6,2);
\draw (4,2) -- (5,1.5);
\draw (3,1.5) -- (4,1);
\draw (2,1) -- (3,0.5);
\draw (7,2.5) node[below right]{$q_1$};
\draw (6,2) node[below right]{$q_{01}$};
\draw (5,1.5) node[below right]{$q_{001}$};
\draw (4,1) node[below right]{$q_{0001}$};
\draw (3,0.5) node[below right]{$q_{00001}$};
\end{tikzpicture}}
\caption{\label{Fig.arbre}Infinite simple comb probabilized context tree.}
\end{figure}
%%%%%%%%%%%%%%%%%%%%%%
It is proved in \cite{ccpp} that in the irreducible case \emph{i.e.}\@ when $q_{0^{\infty}}(0)\not= 1$, $(U_n)_{n\ge 0}$ has a unique stationary probability measure $\pi$ on the set of left-infinite words $\mathcal{L}$ if and only if  $\Theta_1$ is finite. Similarly, if $\Theta_1<\infty$, Proposition \ref{prop:stat-gen} implies that $(X_n,M_n)_{n\in\g Z}$ has a unique invariant probability measure. The following theorem enlights the links between the VLMC $(U_n)_{n\geq 0}$ and the chain $(X_n,M_n)_{n\in\g Z}$ and their respective stationary probability measure.
%%%%%%%%%%%%%%%%
\begin{thm}[\textbf{infinite comb}]
\label{prop:vlmc} 
\ 
\begin{enumerate}[(i)]
\item Let $(X_n,M_n)_{n\in\mathbb{Z}}$ be a
stationary Markov chain valued in
$\{0,1\}\times\overline{\mathbb{N}}^*$, with transition probabilities \eqref{eq:transition} and
\eqref{eq:extension}, with $p_{1,2}=p_{2,1}=1$. We suppose $\Theta_1<\infty$ (where $\Theta_1$ is defined in (\ref{eq:cond-gene})) and $\forall n\in\overline{\mathbb{N}}^*$,
\begin{equation}\label{cas-simple}
\alpha_{2,n}=\alpha_2.
\end{equation}  
We define for all $n\in\mathbb{N}$,
\begin{equation}
\label{eq:prop:vlmc1}
U_n=\ldots X_{n-2}X_{n-1}X_n.
\end{equation}
Then, $(U_n)_{n\ge 0}$ is a stationary variable
length Markov chain associated with the infinite comb with
\begin{equation}\label{eq:prop:vlmc2}
q_1(0)=\alpha_2, \quad q_{0^n1}(1)=\alpha_{1,n},\quad
q_{0^\infty}(1)=\alpha_{1,\infty}.
\end{equation}
The initial distribution is given by
\(
U_0\overset{(d)}{=} \ldots X_{-2}X_{-1}X_0.
\)
\item Conversely consider a stationary VLMC $(U_n)_{n\geq 0}$ satisfying \eqref{eq:prop:vlmc2}.
For $n\geq 0$, define $X_n$ as the last letter of $U_n$ and $(M_n)_{n\ge 0}$ as in \eqref{eq:def:M}.
Then $\Theta_1<\infty$ (where $\Theta_1$ is defined in (\ref{eq:cond-gene})) and $(X_n,M_n)_{n\ge 0}$ is a stationary Markov chain with transitions
\eqref{eq:transition}, \eqref{eq:extension} and \eqref{eq:prop:vlmc2} and initial data $(X_0,M_0)$.  %: $X_0$ is the last letter of the word $U_0$ and $M_0=\inf\{n\ge 0: X_{-n}\neq X_0\}$.
A stationary Markov chain $(X_n,M_n)_{n\in\mathbb{Z}}$ can therefore be defined using the classical procedure of extension from $\mathbb{N}$ to $\mathbb{Z}$.
\end{enumerate}
\end{thm}
\noindent The following tabular resumes the correspondence between these two models and could be considered as a dictionary (in the case: $a_1=0$ and $a_2=1$).
\begin{center}
\begin{tabular}{|c|c|}
\hline
$(X_n, M_n)_{n \in \mathbb{Z}}$ & $(U_n)_{n \geq 0}$\\
\hline
\hline
$\nu$ & $\pi$\\
\hline
$Q\left((a_i,k),(a_i,k+1)\right)=1-\alpha_{i,k}$ & $q_{a_i^ka_j}(a_i)$ for $j\neq i$\\
\hline
$Q\left((a_i,k),(a_j,1)\right)=\alpha_{i,k}$& $q_{a_i^ka_j}(a_j)$ \\
\hline
\end{tabular}
\end{center}
\pff
\begin{enumerate}[(i)]
\item Due to Definition (\ref{eq:prop:vlmc1}) of the process $(U_n)_{n\geq 0}$, for all $s\in\{0,1\}$, the events $\{U_{n+1}=U_{n}s\}$ and $\{X_{n+1}=s\}$ are equal. Therefore $(U_n)_{n\geq 0}$ is a Markov chain as soon as
\[
\delta_{s,u}:=\mathbb{P}(U_{n+1}=U_{n}s\vert U_n=u)=\mathbb{P}(X_{n+1}=s\vert X_{n}=u_0,\, \ldots,X_{n-k}=u_{-k},\ldots)
\]
only depends on $s\in\{0,1\}$ and $u$, where $u=\ldots u_{-1}u_0\in \{0,1\}^{-\mathbb{N}}$.\\
Suppose first that $u_0=1$. Since $M_n\in\mathbb{N}^*$, $(\alpha_{2,n})_{n\ge 1}$ is constant and
$\lpref(u)=1$, then \eqref{eq:transition} and \eqref{eq:prop:vlmc2} imply that
\[
\delta_{s,u}=(1-\alpha_2)\ind{s=1}+\alpha_2 \ind{s=0}=q_1(s).
\]
Let us now consider the case $u_0=0$. Recall (see Proposition \ref{prop:stat-gen}) that
$M_n\in\mathbb{N}^*$. Consequently, there exists $m\in\mathbb{N}^*$ such that $u=\ldots 1 0^m$. Then $M_n=m$, $\lpref(U_n)=0^m1$ and
\begin{align*}
\mathbb{P}(X_{n+1}=s\vert X_{n}=0,M_n=m,\ldots)
=(1-\alpha_{1,m})\ind{s=0}+\alpha_{1,m} \ind{s=1}=q_{0^m1}(s).
\end{align*}
Next, we prove that $(U_n)_{n\geq 0}$ is stationary. Note that \eqref{eq:prop:vlmc1} yields
$U_n=\psi\left((X_{n-i}\right)_{i\ge 0})$ a.s. where $\psi\left((x_{-n}\right)_{n\ge 0})=\ldots x_{-2}x_{-1}x_0$.
Therefore, for any $\ell\ge 0$,
\[
\mathbb{E}\Big[f\Big(U_{n+1-\ell},\ldots,U_{n+1}\Big)\Big]=\mathbb{E}\Big[f\Big(\psi((X_{n+1-\ell-i})_{i\ge 0}),\ldots,\psi((X_{n+1-i})_{i\ge 0})\Big)\Big].
\]
Since $(X_n,M_n)_{n\in\mathbb{Z}}$ is stationary, then $(X_{m+1-i})_{i\ge 0}\overset{(d)}{=}(X_{m-i})_{i\ge 0}$
for any $m\in\mathbb{Z}$. This implies that $(U_n)_{n\geq 0}$ is stationary.\\
\item Let us now assume that $(U_n)_{n\geq 0}$ is a stationary VLMC. Let $x,x'\in\{0,1\}$, $k,k'\ge 1$, $n\in\mathbb{N}$ and
\[
\delta':=\mathbb{P}(X_{n+1}=x',\,M_{n+1}=k'\vert X_{n}=x,\, M_n=k,\ldots).
\]
Then
\begin{align*}
\delta'&=\ind{k'=1 \atop x\neq x'}\mathbb{P}(U_{n+1}=U_nx'\vert U_n=\ldots x'x^k) +\ind{k'=k+1 \atop x= x'}\mathbb{P}(U_{n+1}=U_nx\vert U_n=\ldots (1-x)x^k)\\
&=\ind{k'=1}\Big[ \ind{x=0 \atop x'=1}\mathbb{P}(U_{n+1}=U_n 1\vert U_n=\ldots 10^k)+\ind{x=1 \atop x'=0} \mathbb{P}( U_{n+1}=U_n 0\vert U_n=\ldots 01^k)  \Big]\\
&+\ind{k'=k+1}\Big[ \ind{x=x'=1}\mathbb{P}(U_{n+1}=U_n 1\vert U_n=\ldots 01^k) + \ind{x=x'=0} \mathbb{P}(U_{n+1}=U_n 0\vert U_n=\ldots 10^k)\Big]\\
&=\ind{k'=1}\Big[\ind{x=0\atop x'=1}q_{0^k1}(1)+\ind{x=1 \atop x'=0}q_1(0)  \Big]+\ind{k'=k+1}\Big[ \ind{x=x'=1} q_1(1)+\ind{x=x'=0} q_{0^k1}(0) \Big].
\end{align*}
Using \eqref{eq:prop:vlmc2} we get
\begin{equation*}
\delta'=\ind{k'=1}\Big[\ind{x=0 \atop x'=1}\alpha_{1,k}+\ind{x=1 \atop x'=0}\alpha_2  \Big]+\ind{k'=k+1}\Big[\ind{x=x'=1}(1-\alpha_2)+\ind{x=x'=0}(1-\alpha_{1,k}) \Big].
\end{equation*}
Then \eqref{eq:transition} follows directly with $\alpha_{2,n}=\alpha_2$. \QED
\end{enumerate}
The following result is a corollary of Proposition~\ref{prop:stat-gen}. It enables us to compare the expression of the invariant measure {\color{red}$\nu$} from the model $(X_n,M_n)_{n\in \mathbb{Z}}$ with the invariant measure $\pi$ for the VLMC \emph{infinite comb} (see Section~\ref{app:double-peigne} in the Appendix for notations about VLMC).
\begin{cor}
 \label{cor:stat-gen}
Under the condition $\Theta_1<\infty$, there exists a unique invariant probability measure $\nu$ for the Markov chain $(X_n,M_n)$ given by $\nu(a_1,\infty)=\nu(a_2,\infty)=0$ and for all $m\geq 1$,
\begin{equation}
\label{inv-mesure-peigne}
\nu(a_1,m)=\frac{1}{\Theta_1+\Theta_2}\mathcal{P}_1(m) \quad \mbox{and}\quad \nu(a_2,m)=\frac{\alpha_2(1-\alpha_2)^{m-1}}{1+\alpha_2\Theta_1},
\end{equation}
where $\Theta_1=1/\alpha_2$.
In particular one gets
\[\nu(a_2,\mathbb{N}^*)=\pi(a_2)=\frac{1}{1+\alpha_2\Theta_1}.\]
\end{cor}
\pff 
Proposition~\ref{prop:stat-gen} with 
\[\mathcal{P}=\left( \begin{array}{cc} 0 & 1\\ 1 & 0\end{array}  \right)\quad\mbox{and}\ v^*=\frac{1}{2}\,(1,1)\]
lead to \eqref{inv-mesure-peigne} and 
\[\nu(a_2,\mathbb{N}^*)=\sum_{m\ge 1}\nu(a_2,m)=\frac{1}{1+\alpha_2\Theta_1}.\]
Consequently one has
\[\nu(a_2,\mathbb{N}^*)=\frac{1}{1+q_1(0)\sum_{n\geq 1}\prod_{k=1}^{n-1}(1-q_{0^k1}(1))}=\frac{1}{\sum_{n\geq 0}\prod_{k=0}^{n-1}q_{0^k1}(0)}=\pi(a_2), \]
which is fortunately (!) the invariant measure obtained in \cite{ccpp}. \QED

%%%%%%%%%%%%%%%%%%%%%%%%%%%%%%%%%%%%%%%%%%%%%%%%%%%%%%%%%%%%%%%
%
%       DOUBLE COMB
%
%%%%%%%%%%%%%%%%%%%%%%%%%%%%%%%%%%%%%%%%%%%%%%%%%%%%%%%%%%%%%%%
\subsection{The double infinite comb}
\begin{figure}[htbp]\label{fig.arbredouble}
\centerline{\begin{tikzpicture}[scale=1]
\draw (2,1) -- (5,2.5) ;
\draw (5,2.5) -- (8,1);
%\draw [>=stealth,->] (0,0) -- (0,1) ;
\draw [dashed] (0,0) -- (2,1) ;
\draw [dashed] (8,1) -- (10,0);
\draw (4,2) -- (4.75,1.25);
\draw (3,1.5) -- (3.75,0.75);
\draw (2,1) -- (2.75,0.25);
\draw (1,0.5) -- (1.75,-0.25);
\draw (6,2) -- (5.25,1.25);
\draw (7,1.5) -- (6.25,0.75);
\draw (8,1) -- (7.25,0.25);
\draw (9,0.5) -- (8.28,-0.25);
%\draw (4,2) node[above left]{$0$};
%\draw (3,1.5) node[above left]{$00$};
%\draw (2,1) node[above left]{$000$};
%\draw (1,0.5) node[above left]{$0000$};
\draw (4.75,1.25) node[below]{$q_{01}$};
\draw (3.75,0.75) node[below]{$q_{001}$};
\draw (2.75,0.25) node[below]{$q_{0001}$};
\draw (1.75,-0.25) node[below]{$q_{0^n1}$};
\draw (0,0) node[below]{$q_{0^\infty}$};
%\draw (6,2) node[above right]{$1$};
%\draw (7,1.5) node[above right]{$11$};
%\draw (8,1) node[above right]{$111$};
%\draw (9,0.5) node[above right]{$1111$};
\draw (5.25,1.25) node[below]{$q_{10}$};
\draw (6.25,0.75) node[below]{$q_{110}$};
\draw (7.25,0.25) node[below]{$q_{1110}$};
\draw (8.25,-0.25) node[below]{$q_{1^n0}$};
\draw (10,0) node[below]{$q_{1^\infty 0}$};
\end{tikzpicture}}
\caption{\label{figPeigneInfini}
infinite double comb probabilized context tree.}
\end{figure}
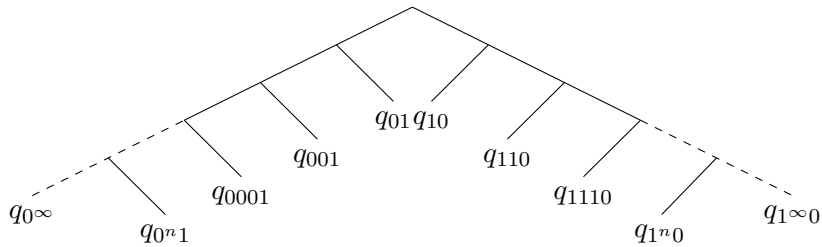

Let us now present the double
infinite comb. 
%%%%%%%%%%%%%%%%%%
Consider the probabilized context tree given on Figure~\ref{figPeigneInfini} (hereafter called double infinite comb).
In this case, there are two infinite leaves $0^\infty$ and $1^{\infty}$ and countably many finite leaves $0^n1$ and $1^n0$, $n\in\mathbb{N}$, so that
\[
\mathcal{C}=\{0^n1,\, n\ge 1\}\cup \{ 1^n0,\, n\ge 1 \}\cup \{0^\infty\}\cup \{1^\infty\}.
\]
The data of a corresponding VLMC consists thus in Bernoulli probability measures on $\{ 0,1\}$:
\[q_{0^\infty}, q_{1^\infty}, 
{\rm ~and~~} q_{0^n1}, q_{1^n0},~n\in \mathbb{N}^*.\]
We refer to Appendix~\ref{app:double-peigne} to see that the finiteness of $\Theta_1$ and $\Theta_2$ implies the existence of a unique invariant measure for this VLMC. 
\begin{thm}[\textbf{double infinite comb}]
\label{prop:vlmcd}
\begin{enumerate}[(i)]
\item Let $(X_n,M_n)_{n\in\mathbb{Z}}$ be a stationary Markov chain with
transition probabilities \eqref{eq:transition} and
\eqref{eq:extension}. We suppose $\Theta_1<\infty$ and $\Theta_2<\infty$ (where $\Theta_i$ is defined in (\ref{eq:cond-gene})). Then the process $(U_n)_{n\geq 0}$ defined by
\eqref{eq:prop:vlmc1} is a stationary variable length Markov chain
associated with the double infinite comb with
\begin{equation}\label{eq:prop:vlmc2d}
q_{1^n0}(0)=\alpha_{2,n},\quad q_{1^\infty}(0)=\alpha_{2,\infty},
\quad q_{0^n1}(1)=\alpha_{1,n},\quad
q_{0^\infty}(1)=\alpha_{1,\infty}.
\end{equation}
The initial data is given by $U_0=\ldots X_{-2}X_{-1}X_0$.\\
\item Conversely let $(U_n)_{n\ge 0}$ be a stationary VLMC satisfying \eqref{eq:prop:vlmc2d}. For $n\geq 0$, define $X_n$ by the last letter of $U_n$ and $(M_n)_{n\ge 0}$ as in \eqref{eq:def:M}.
Then $(X_n,M_n)_{n\ge 0}$ is a stationary Markov chain with transitions \eqref{eq:transition}, \eqref{eq:extension} and \eqref{eq:prop:vlmc2d} and with initial data $(X_0,M_0)$. This stationary Markov chain can be extended on the time space $\mathbb{Z}$ as usual.
\end{enumerate}
\end{thm}
\noindent The arguments for the proof are similar to those presented in Theorem \ref{prop:vlmc}.

\medskip

\noindent As for the simple infinite comb, the invariant measure of the first margin of the Markov chain $(X_n,M_n)_{n\in\mathbb{Z}}$ corresponding to the double infinite comb can be compared with the invariant measure $\pi$ for the VLMC \emph{double infinite comb} calculated in Appendix~\ref{app:double-peigne}.
\begin{cor}
 \label{cor:inv-mes-double}
Under the condition $\Theta_1<\infty$ and $\Theta_2<\infty$, there exists a unique invariant probability measure $\nu$ for the Markov chain $(X_n,M_n)$ given by $\nu(a_1,\infty)=\nu(a_2,\infty)=0$ and for all $m\geq 1$ and $i=1,2$,
\begin{equation}
\label{inv-mes-double}
 \nu(a_i,m)=\frac{1}{\Theta_1+\Theta_2}\mathcal{P}_i(m). \end{equation}
Consequently one gets
\[\nu(a_2,\mathbb{N}^*)=\frac{\Theta_2}{\Theta_1+\Theta_2}=\pi(a_2).\]
\end{cor}
\pff
Again \eqref{inv-mes-double} is a direct consequence of Proposition~\ref{prop:stat-gen}. Suming up it comes
\[\nu(a_2,\mathbb{N}^*)=\frac{\Theta_2}{\Theta_1+\Theta_2},\]
with 
\[\Theta_1=\sum_{n\geq 1}\prod_{k=1}^{n-1}(1-q_{0^k1}(1))=\sum_{n\geq 0}\prod_{k=1}^{n}q_{0^k1}(0)\]
and
\[\Theta_2=\sum_{n\geq 1}\prod_{k=1}^{n-1}(1-q_{1^k0}(0))=\sum_{n\geq 0}\prod_{k=1}^{n}q_{1^k0}(1),\]
which is exactly the calculation of $\pi(a_2)$ in Appendix~\ref{app:double-peigne}.
\begin{rem}
 \label{rem:gene:alphabet}
The results developed in Theorem \ref{prop:vlmc} and Theorem \ref{prop:vlmcd} can be generalized to context trees
which are based on a finite alphabet $\{a_1,\ldots,a_K\}$ and composed with a finite number of combs.
The corresponding Markov chain $(X_n,M_n)_{n\in\g Z}$ is then valued in the state space
$\{a_1,\ldots,a_K\}\times\overline{\mathbb{N}}^*$.
\end{rem}
Of particular interest are variable length Markov chains $(U_n)$ associated with the infinite comb or the double infinite comb. While the sequence $(X_n)$ formed by the last letters of the process $U_n= \ldots X_{n-1}X_n$ is not a Markov process, except for very particular $q_c$, % in the general case. 
the previous theorems show that it suffices to add a memory process $(M_n)$ to get a Markov chain $(X_n,M_n)$. Note that $(U_n)$ takes its value in the non-countable space $\mathcal{L}$ and Theorems~\ref{prop:vlmc} and \ref{prop:vlmcd} allow to associate by a one to one correspondence a Markov chain $(X_n,M_n)$ which is valued in the countable set $\{0,1\}\times \overline{\mathbb{N}}^*$. This reduction of the size of the state space (which becomes here minimal) is made possible by the particular shape of the context tree: for instance, the VLMC associated with the bamboo blossom defined in \cite{ccpp} is not equivalent to a Markov Chain $(X_n, M_n)$ with a real memory process. Nevertheless for suitable VLMC we suggest to introduce the following application 
\[(U_n)_{n}\mapsto \left(\lpref(U_n)\right)_{n},\]
which should permit to generalize the reduction of the state space.
The image process is not Markovian in the general case, even under the stationary distribution for $U_n$. A conjecture: the process $\left(\lpref(U_n)\right)_n$ is Markovian (and thus defines an automaton) if and only if the associated context tree has a completeness property, studied in a companion paper \cite{ccpp3}.
%%%%%%%%%%%%%%%%%%%%%%%%%%%%%%%%%%%%%%%%%%%%%%%%%%%%%
%%
%%          Paths description-Distribution
%%
%%%%%%%%%%%%%%%%%%%%%%%%%%%%%%%%%%%%%%%%%%%%%%%%%%%%%
\section{Distribution of the persistent random walk}
\label{sect:persistent}
\setcounter{equation}{0}
By definition, a random walk $(S_n)_{n\ge 0}$ is a process whose increments are independent. It is often pertinent, for instance in modeling, to begin with the increments and second to study the associated random walk. Let us give an example coming from finance. Suppose that $S_n$ is the price at time $n$ of an asset. In the Cox, Ross and Rubinstein model, the non-arbitrage condition implies that the relative increments $\left(\frac{S_n-S_{n-1}}{S_{n-1}};\ n\ge 1 \right)$ are independent.\\
We study here a class of \emph{additive} processes $(S_n)$ of the type
\begin{equation}
\label{eq:def-persist-part}
S_n=\sum_{k=0}^n X_k,\quad n\ge 0,
 \end{equation}
 where the increments $(X_n)$ are not independent. A tentative of considering increments with short dependency has been already developed in \cite{Tapiero-Vallois} and \cite{Tapiero-Vallois2}. In these studies, the authors have supposed that $(X_n)$ is a Markov chain. We would like to go further here introducing \emph{variable length memory} between the increments.\\
 We consider in this section, a Markov chain $(X_n,M_n)_{n\ge 0}$ with transition probability \eqref{eq:transition} and
\eqref{eq:extension} and we assume 
\[
K=2,\, a_1=-1,\, a_2=1. 
\]
The process $(S_n)$ defined by \eqref{eq:def-persist-part} is called a \emph{persistent random walk}. This terminology comes from \cite{eckstein00}.\\
%Equivalently, consider a stationary VLMC $(U_n)_{n\geq 0}$ associated with a double infinite comb and its letter process $(X_n)_{n\geq 0}$.
%Let us introduce the persistent random walk 
%\begin{equation*}
%%\label{eq:def-persist-part}
%S_n=\sum_{k=0}^n X_k,\quad n\ge 0.
% \end{equation*}
% Classical random walks are crucial processes in many domains of application: they represents the simplest random movement. The trajectory consists of taking successive independent random steps. In the financial framework, it is sometimes possible to represent the price of a fluctuating stock by a random walk. However we have to assume that all the financial agents act in a markovian way: they are memoryless. This point of view is perhaps simplistic! Considering a random walk whose increments are not independent but consist in a one-order Markov chain permits to introduce a short memory (in financial models for instance) and was the subject of \cite{Herrmann-Vallois}. In this section, we focus our attention to increments which corresponds to the last letter of a VLMC. In other words we introduce long memory in the random walk. We shall describe the distribution of this new random process and in particular observe if the memory have a important influence in the large time behaviour of the walk.

A path description of $(S_n)$ is given in Section \ref{subsec:paths1}, puting ahead the breaking times $(T_n)_{n\geq 1}$. We give in Section \ref{sect:distr} the explicit distribution of $S_n$. Although the law of $S_n$ is complicated, we can determine explicitely the distribution and the generating function of the position of the persistent random walk at an exponential independent random time. The double generating funtion will play an important role in Section~\ref{sec:GITN}. 

We end this section studying how $S_n$ fluctuates as $n\to\infty$. Indeed it is not so far from the persistent walk with one-order Markovian increments. We prove a law of large number and a central limit theorem. We recover the classical setting where $(X_n)$ is a Markov chain. We have introduced variable memory to $(X_n)$, but it seems that it is not sufficient to obtain new asymptotic behavior: it would be therefore very interesting to investigate the behaviour of the random walk when \emph{mixing assumptions} are relaxed, \emph{i.e.} when the length of the memory increases significantly to give a real \emph{persistent memory effect to the random walks}.  

%%%%%%%%%% 
%{\color{pink}B:  ici ou plus tard section 4.2. : id\'ee du d\'eroulement de cette section. 
%
%SI JE COMPRENDS BIEN : 
%
%dans 4.2 : on calcule dans 4.2.1 la loi de $S_n$ dans Prop 4.1 puis dans 4.2.2 la loi de $S_{\tau}$ dans theorem 4.10. Donc aussi celle de $N_n(1)$ puisque c'est pareil que $S_n$ grace \`a \eqref{eq:prop:loicm1}, et celle de $N_{\tau}$ puisque c'est pareil que $S_{\tau}$. On connait plus ou moins explicitement aussi la loi de $T_{2m}$. On comprend au fin fond des remarques 2.8 et 4.6 que la m\'emoire $M_n$ c'est en gros $T_{n+1}-T_n$ et sa loi est donn\'ee par les $\rond P_i(n)$.
%
%Puis dans la section 4.2.3 et avec Prop 4.14, on a une asymptotique de $S_n$.
%
%Puis dans la section 4.3, on definit the time continuous process $S^{\eps}(t)$ et on montre dans le theorem 4.17 qu'il converge en montrant dans la Prop 4.16 que son counting process $N_t^{\eps}$ converge. Et parce qu'on a un lien entre $S_n$ et son counting process grace \`a \eqref{eq:rel-S-N} et le m\^eme lien en temps continu.
%
%Question : est-ce que la section 4.3 utilise la section 4.2 ? 
%}
%%%%%%%%%

\subsection{Paths description}
\label{subsec:paths1} Since $X_n$ is $\{-1,1\}$-valued, it is
clear that the trajectory of $(S_n)_{n \geq 0}$ is a sequence of straight
lines with slopes $\pm 1$, and the instants of breaks are
$(T_n)_{n\ge 1}$ which were introduced in \eqref{eq:changt}. 

\medskip
\noindent Let us assume that $S_0=X_0=1$, then the trajectory increases step $1$ by step $1$ till $T_1-1$ where it reaches a first local maximum. After that time, it decreases and reaches a local minimum at time $T_2-1$ and so on. The trajectory of $(S_n)_{n\in\mathbb{N}}$ corresponds to the linear interpolation between
the sequence of points $(W_n,Z_n)_{n\ge 0}$ where $W_0=0$, $Z_0=1$ and for $n\ge 1$,
\[(W_n,Z_n)=\left(T_n-1, S_{T_n}-(-1)^n\right)=\left(T_n-1, \sum_{k= 1}^n(-1)^{k-1}(T_k-T_{k-1})\right).\]
%%%%%%%%%%%%%%%%%%%%%%%%%%%%%%%%%%%%%%%%%%%%%%%%%%%%%%%%%%%%%%%%%%%%%%%%%%%%%%%%%%%%%%%%%%%%%%%%%%%%%%%
\begin{figure}[!h]
\begin{minipage}{0.52\linewidth}
\centerline{\begin{tikzpicture}[scale=0.9]
\draw (0,-2) -- (0,3) ;
\draw [>=stealth,->] (0,0) -- (8,0) ;
\draw (0,0) -- (8,0) ;
\draw [dotted, very thin, gray] (0,-2) grid (8,3);
\draw (8,0) node[right]{$n$} ;
\draw (0,3) node[left]{$S_n$};
\draw [>=stealth,->] (0,1) -- (1,2) ;
\draw [>=stealth,->] (1,2) -- (2,3) ;
\draw [>=stealth,->] (2,3) -- (3,2) ;
\draw[dashed] (3,2) -- (3,0) ;
\draw (3, 0) node[below]{$T_{1}$};
\draw [>=stealth,->] (3,2) -- (4,1) ;
\draw [>=stealth,->] (4,1) -- (5,0) ;
\draw [>=stealth,->] (5,0) -- (6,-1) ;
\draw [>=stealth,->] (6,-1) -- (7,0) ;
\draw[dashed] (7,0) -- (7,0) ;
\draw (7, 0) node[above]{$T_{2}$};
\draw [>=stealth,->] (7,0) -- (8,1) ;
\end{tikzpicture}}
\centerline{$X_1=\ldots=X_8=11(-1)(-1)(-1)(-1)11$}
\end{minipage}
\begin{minipage}{0.52\linewidth}
\centerline{\begin{tikzpicture}[scale=0.9]
\draw (0,-2) -- (0,3) ;
\draw [>=stealth,->] (0,0) -- (7,0) ;
\draw (0,0) -- (7,0) ;
\draw [dotted, very thin, gray] (0,-2) grid (7,3);
\draw (7,0) node[right]{$n$} ;
\draw (0,3) node[left]{$S_n$};
\draw [>=stealth,->] (0,-1) -- (1,-2) ;
\draw [>=stealth,->] (1,-2) -- (2,-1) ;
\draw[dashed] (2,-1) -- (2,0) ;
\draw (2, 0) node[above]{$T_{1}$};
\draw [>=stealth,->] (2,-1) -- (3,0) ;
\draw [>=stealth,->] (3,0) -- (4,1) ;
\draw [>=stealth,->] (4,1) -- (5,2) ;
\draw [>=stealth,->] (5,2) -- (6,1) ;
\draw[dashed] (6,0) -- (6,1) ;
\draw (6, 0) node[below]{$T_{2}$};
\draw [>=stealth,->] (6,1) -- (7,0) ;
\end{tikzpicture}}
\centerline{$X_1\ldots X_7=(-1)1111(-1)(-1)$,}
\end{minipage}

\caption{\label{Fig1}Trajectories of $(S_t)_{t\geq 0}$ when either $S_0=1$ or $S_0=-1$.}
\end{figure}
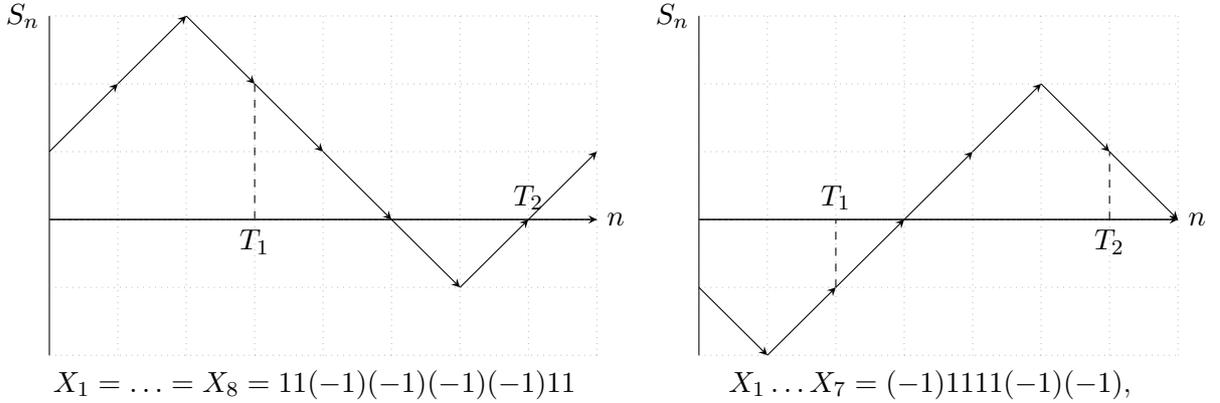
%%%%%%%%%%%%%%%%%%%%%%%%%%%%%%%%%%%%%%%%%%%%%%%%%%%%%%%%%%%%%%%%%%%%%%%%%%%%%%%%%%%%%%%%%%%%%%%%%%%%%%
If $S_0=X_0=-1$, then the behaviour of the process is similar
and reduces on a succession of increasing and decreasing parts. The trajectory 
$(S_t)_{t\ge 0}$ is a linear interpolation between
$(W_n,Z_n')_{n\ge 0}$ where $W_0=0$, $Z_0'=-1$ and for $n\ge 1$,
\[(W_n,Z_n')=\left(T_n-1, S_{T_n}+(-1)^{n}\right)=\left(T_n-1, \sum_{k=
1}^n(-1)^{k}(T_k-T_{k-1})\right).\]
Note that $Z'_n=-Z_n$.

\medskip

\noindent Let us introduce the counting process $(N_t)_{t\in\mathbb{N}}$ whose jump times are $T_n$:
\begin{equation}
\label{eq:def:compt}
N_t=\sup\{ n\ge 1:\ T_n\le t \} =\sum_{n\ge 1}\ind{T_n\le t},\quad t\in\mathbb{N} .
\end{equation}
%Obviously, \eqref{eq:def:compt} is equivalent to
%\begin{equation}\label{eq:equiv-compt}
%N_t=\sum_{n\ge 1}1_{\{ T_n\le t \}},\quad t\in\mathbb{N}.
%\end{equation}
From now on, we suppose that $(X_0,M_0)=(1,1)$. Note that the case $(X_0,M_0)=(-1,1)$ can be deduced from the former case changing
$X$ in $-X$. 

\noindent The counting process $(N_t)_{t\geq 0}$ will play an important role in the study of $(S_t)_{t\geq 0}$ (see  Section~\ref{sec:GITN}) and $(S_t)_{t \geq 0}$ can be expressed via $(N_t)_{t \geq 0}$ as:
\begin{equation}
\label{eq:rel-S-N}
S_t=\sum_{n=0}^{t}(-1)^{N_n}.
\end{equation}
There is a one-to-one correspondence between $(M_t)_{t \geq 0}$
and $(T_n)_{n\ge 0}$:
\begin{equation}\label{eq:oneto}
\{k;\ M_k=1\}=\{T_n;\ n\ge 0\}.
\end{equation}
$\{N_s;\ s\le t\}$ can be expressed via $\{M_s;\ s\le t\}$ and vice
and versa. Indeed, \eqref{eq:def:compt} \eqref{eq:oneto} and
\eqref{eq:def:M} imply
\begin{equation}\label{eq:N-M}
N_t=\sum_{k=1}^t\ind{M_k=1}\quad\mbox{and}\quad M_t=1+\sup\{n\ge
0:\ N_{t-n}=N_t\},\, t\in\mathbb{N}.
\end{equation}
%%%%%%%%%%%%%%%%%%%%%%%%%%%%%%%%%%%%%%%%%%%%%%%%%%%%%%%%%%%%%%%%
%%
%%
%%%%%%%%%%%%%%%%%%%%%%%%%%%%%%%%%%%%%%%%%%%%%%%%%%%%%%%%%%%%%%%%
%\subsection{Distribution of a persistent random walk}
\subsection{Distribution of the persistent random walk at a fixed time}
\label{sect:distr}
In this section we give the explicit distribution of the persistent random walk at any fixed time.\\
We recall that $(X_n, M_n)_{n\ge 0}$ is a $\{-1,1\}\times\overline{\mathbb{N}}^*$-valued Markov chain with transitions matrix $Q$ defined by \eqref{eq:transition} and starting values $(X_0,M_0)=(1,1)$. Therefore the law of $(X_n,M_n)$ is given by $Q^n$. However the calculation of $Q^n$ is untractable. This leads to restrict ourselves to the law of $X_n$.\\
Let us define
\[
\mathcal{N}(m,b):=\Big\{u\in(\mathbb{N}^*)^m:\ u_1+\ldots+u_m=b\Big\},\quad m\ge 1,\ b\ge 1.
\]
and 
\begin{equation}
\label{eq:defdeA}
A_i(m,b)=\sum_{u\in \mathcal{N}(m,b)}\mathcal{P}_i(u_1)\ldots\mathcal{P}_i(u_m)\ \alpha_{i,u_1}\times\ldots\times \alpha_{i,u_m},\quad i=1,2
\end{equation}
with $A_i(m,b)=0$ for $0\le b<m$ and $A_i(0,b)=\ind{b=0}$.\\
The distribution of the random walk $(S_n)_{n \geq 1}$ can be directly linked to the occupation measure $L_n(1)$ of the increments $(X_n)_{n \geq 1}$ in the following way:

\begin{prop} {\mathversion{bold} \bf (distribution of $S_n$)\mathversion{normal}}
\label{prop:loicm}
Suppose that $(X_0,M_0)=(1,1)$.
\begin{enumerate}[(i)]
\item Let us introduce the local time
\begin{equation}
\label{eq:prop:loicm}
L_n(1):=\sum_{k=1}^n\ind{X_k=1},
\end{equation}
then the random walk satisfies for $n\ge 1$,
\begin{equation}
\label{eq:prop:loicm1}
S_n=1+2L_n(1)-n.
\end{equation}
Consequently, for any $0\le k\le n$:
\begin{equation}\label{eq:prop:loicm21}
\eta_n(k):=\mathbb{P}(S_n=1+2k-n)=\mathbb{P}(L_n(1)=k).
\end{equation}
\item Moreover, for $0\le k\le n$, we have 
\begin{equation}
\label{eq:distrSn}
\eta_n(k)=\eta_n^{(1)}(k)+\eta_n^{(2)}(k)
\end{equation}
with
\begin{align}
\label{eq:prop:align}
\eta_n^{(1)}(k)&=\sum_{1\le m\le (k+1)\wedge (n-k) } A_2(m,k+1)\sum_{\ell=1}^{n-k-m+1}A_1(m-1,n-k-\ell)\mathcal{P}_1(\ell)\\
\label{bis}
\eta_n^{(2)}(k)&=\sum_{0\le m\le k\wedge (n-k) } A_1(m,n-k)\sum_{\ell=1}^{k-m+1}A_2(m,k+1-\ell)\mathcal{P}_2(\ell).
\end{align}
\end{enumerate}
\end{prop}
\pff
\begin{enumerate}[(i)]
\item Using the definition of $L_n(1)$, it comes
\begin{align*}
S_n&=1+\sum_{i=1}^n\ind{X_i=1}-\sum_{i=1}^n\ind{X_i=-1}\\
&=1+\sum_{i=1}^n\ind{X_i=1}-\left(n-\sum_{i=1}^n\ind{X_i=1}\right)\\
&=1+2L_n(1)-n.
\end{align*}
\item In order to compute $\eta_n(k)$, it is convenient to use the family of stopping times $(T_n)$ introduced in \eqref{eq:changt}. The probability of the event $\{L_n(1)=k\}$ can be decomposed into two parts, according to the fact that time $n$ arrives on the way up or on the way down: 
\begin{equation}\label{eq:proof}
\eta_n(k)=\sum_{m\ge 0}\eta_n^{(1)}(k,m)+\sum_{m\ge 1}\eta_n^{(2)}(k,m)
\end{equation}
where
\[
\eta_n^{(1)}(k,m):=\mathbb{P}\Big(L_n(1)=k, \, T_{2m}\le n<T_{2m+1}\Big),\quad m\ge 0
\]
and
\[
\eta_n^{(2)}(k,m):=\mathbb{P}\Big(L_n(1)=k, \, T_{2m-1}\le n<T_{2m}\Big),\quad m\ge 1.
\]
\paragraph{\mathversion{bold}First step --- Computation of $\eta_n^{(1)}(k,m)$ for $n\ge k$.\mathversion{normal}}

Suppose first that $m\ge 1$. On the set $\{L_n(1)=k,\, T_{2m}\le n<T_{2m+1}\}$, we define for $0\le i<m$, the length of the $i$\textsuperscript{th} ascent $W_i:=T_{2i+1}-T_{2i}$, $W_m:=n+1-T_{2m}$ and the length of the $i$\textsuperscript{th} descent  $V_i:=T_{2i}-T_{2i-1}$ for $1\le i\le m$. Then
\begin{equation}\label{eq:decosom}
W_0+W_1+\dots+W_m+V_1+\ldots+V_m=n+1,\quad W_0+W_1+\ldots+W_m=k+1.
\end{equation}
Therefore for $\underline{W}=(W_0,\ldots,W_m)$ and $\underline{V}=(V_1,\ldots,V_m)$ we get
\begin{align}
\label{eq:produit}
\eta_n^{(1)}(k,m)&=\sum_{w\in\mathcal{N}(m+1,k+1)}\quad\sum_{v\in\mathcal{N}(m,n-k)}\mathbb{P}(\underline{W}=w,\, \underline{V}=v).
\end{align}
Using the distributions of $T_{2i+1}-T_{2i}$ and $T_{2i+2}-T_{2i+1}$ given in Proposition \ref{prop:changtpente-sym}, we obtain
\begin{align}\label{eq:calc}
\mathbb{P}(\underline{W}=w,\, \underline{V}=v)&=\mathcal{P}_2(w_1)\alpha_{2,w_1}\mathcal{P}_1(v_1)\alpha_{1,v_1}\times\ldots\nonumber\\
&\times \mathcal{P}_2(w_m)\alpha_{2,w_m}\mathcal{P}_1(v_m)\alpha_{1,v_m}\mathcal{P}_2(u_{m+1}).
\end{align}
It is clear that \eqref{eq:produit} and \eqref{eq:calc} imply
\begin{align}
\label{eq:etap11}
\eta_n^{(1)}(k,m)=\widehat{A}_2(m+1,k+1) A_1(m,n-k),
\end{align}
where $A_1$ is defined by \eqref{eq:defdeA} and for $m\ge 2$ and $i\in \{1,2\}$,
\begin{equation}
\label{eq:defdehatA}
\widehat{A}_i(m,b):=\sum_{w\in\mathcal{N}(m,b)}\mathcal{P}_i(w_1)\times\ldots\times\mathcal{P}_i(w_m)\ \alpha_{i,w_1}\times\ldots\times \alpha_{i,w_{m-1}},
\end{equation}
and $\widehat{A}_i(1,b)=\mathcal{P}_i(b)$. \\ If $m=0$, then $n=k$, $\eta^{(1)}_n(k,m)=\mathcal{P}_2(n+1)$. Therefore \eqref{eq:etap11} holds with $m=0$.\\

\paragraph{\mathversion{bold} Step 2 --- Computation of $\eta_n^{(2)}(k,m)$.\mathversion{normal}}

Similarly, define on $\{ L_n(1)=k,\, T_{2m-1}\le n<T_{2m} \}$, $W_i:=T_{2i+1}-T_{2i}$ for $0\le i<m$, $V_i:=T_{2i}-T_{2i-1}$ for $1\le i< m$ and $V_m:=n+1-T_{2m-1}$ then:
\[
W_0+W_1+\dots+W_{m-1}+V_1+\ldots+V_m=n+1,\quad W_0+W_1+\ldots+W_{m-1}=k+1.
\]
For $\underline{W}=(W_0,\ldots,W_{m-1})$ and $\underline{V}=(V_1,\ldots,V_m)$ we get
\begin{align}
\label{eq:produit2}
\eta_n^{(2)}(k,m)&=\sum_{w\in\mathcal{N}(m,k+1)}\quad\sum_{v\in\mathcal{N}(m,n-k)}\mathbb{P}(\underline{W}=w,\, \underline{V}=v)\nonumber\\
&=A_2(m,k+1)\widehat{A}_1(m,n-k).
\end{align}
Combining \eqref{eq:proof}, \eqref{eq:etap11} and \eqref{eq:produit2} leads to
\begin{equation}\label{eq:fini}
\eta_n(k)=\sum_{m=0}^{k\wedge(n-k)}\widehat{A}_2(m+1,k+1)A_1(m,n-k)+\sum_{m=1}^{(k+1)\wedge(n-k)}A_2(m,k+1)\widehat{A}_1(m,n-k)\end{equation}
In order to prove \eqref{eq:distrSn}, it suffices to express $\widehat{A}_i$ in terms of $A_i$. For $b\ge m> 1$ we observe that
\[
\mathcal{N}(m,b)=\Big\{(w,w_m):\ w \in\mathcal{N}(m-1,j),\, w_m=b-j,\, m-1\le j\le b-1\Big\}.
\]
Hence, for $b\ge m>1$,
\begin{equation}\label{eq:achapa}
\widehat{A}_i(m,b)=\sum_{j=m-1}^{b-1}A_i(m-1,j)\mathcal{P}_i(b-j)=\sum_{\ell=1}^{b-m+1}A_i(m-1,b-\ell)\mathcal{P}_i(\ell).
\end{equation}
Observe that \eqref{eq:achapa} is still valid if $m=1$, since $A_i(0,b)=\ind{b=0}$ and $\widehat{A}_i(1,b)=\mathcal{P}_i(b)$.

The decomposition \eqref{eq:achapa} permits to transform \eqref{eq:fini} into \eqref{eq:distrSn}. \QED
\end{enumerate}
\begin{rem}\label{cor:peigsim}
In the particular situation $\alpha_{2,k}=\alpha_2$ for any $k\ge 1$, which is associated with the simple infinite comb (Section~\ref{sec:contexttree}), then the distribution of $S_n$ given by \eqref{eq:prop:loicm21} and \eqref{eq:distrSn} can be simplified since
\[A_2(m,b)=\binom{b-1}{m-1}(1-\alpha_2)^{b-m}\alpha_2^m=\alpha_2\widehat{A}_2(m,b),\quad b\ge m\ge 1.\]
\end{rem}
Of course by symmetry we get also a similar expression of $A_1$ if $\alpha_{1,k}=\alpha_1$ for any $k\ge 1$. Combining both identities, Proposition~\ref{prop:loicm} gives the distribution of $S_n$ when $X_n$ is a Markov chain. Let us just note that the associated VLMC is very particular and the generating function of $S_n$ was already presented in \cite{Tapiero-Vallois}.
%\pff
%Since $\alpha_{2,k}=\alpha_2$ for $k\ge 1$, we obtain $\mathcal{P}_2(k)=(1-\alpha_2)^{k-1}$ and
%\begin{align*}
%A_2(m,b)&=\sum_{u\in\mathcal{N}(m,b)}\mathcal{P}_2(u_1)\ldots\mathcal{P}_2(u_m)\ \alpha_{2,u_1}\times\ldots\times\alpha_{2,u_m}\\
%&=\sum_{u\in\mathcal{N}(m,b)}(1-\alpha_2)^{u_1+\ldots +u_m-m}\alpha_2^m\\
%&=(1-\alpha_2)^{b-m}\alpha_2^m\Big(\# \mathcal{N}(m,b)
%\Big)
%\end{align*}
%\bnote{notation $|\ |$}Using
%\[
%\sum_{b\ge m} x^{b-m}\Big(\#\mathcal{N}(m,b)
%\Big)=\sum_{b\ge m}\sum_{u_1,\ldots,u_m}x^{u_1+\ldots+u_m-m}\ind{u_1+u_2+\ldots+u_m=b}=\left(\frac{1}{1-x}\right)^m,
%\]
%we deduce easily that $\#\mathcal{N}(m,b)=\binom{b-1}{m-1}$. \bnote{aussi r\'esultat classique combi}Therefore
%\begin{align*}
%A_2(m,b)=\binom{b-1}{m-1}(1-\alpha_2)^{b-m}\alpha_2^m.
%\end{align*}
%Furthermore
%\[\widehat{A}_2(m,b)=\sum_{u\in\mathcal{N}(m,b)}\mathcal{P}_2(u_1)\ldots\mathcal{P}_2(u_m)\ \alpha_{2}^{m-1}=\frac{A_2(m,b)}{\alpha_2}. \quad \blacksquare\]
%In \cite{Tapiero-Vallois}, the authors have determined the generating function of $S_n$ when $X_n$ is a Markov chain. {\color{red}Although this case just corresponds to a particular VLMC,} it is interesting to note that Proposition \ref{prop:loicm} gives the distribution of $S_n$.
\begin{cor}\label{cor:peigtriv} Suppose that $\alpha_{1,k}=\alpha_1$ and $\alpha_{2,k}=\alpha_2$ for any $k\ge 1$. This means that $(X_n)$ is a $\{ -1,1 \}$-valued Markov chain with transition matrix 
\begin{align*}
\left(\begin{array}{cc}
1-\alpha_1 & \alpha_1\\
\alpha_2 & 1-\alpha_2
\end{array}
\right).
\end{align*}
Then one has
\begin{align*}
\mathbb{P}(L_n(1)=k)&=\sum_{m=1}^{(k+1)\wedge (n-k)}\binom{k}{m-1}\binom{n-k-1}{m-1}\alpha_1^{m-1}(1-\alpha_{1})^{n-k-m}\alpha_2^m(1-\alpha_2)^{k+1-m}\\
&+\sum_{m=1}^{k\wedge (n-k)}\binom{k}{m}\binom{n-k-1}{m-1}\alpha_1^{m}(1-\alpha_{1})^{n-k-m}\alpha_2^m(1-\alpha_2)^{k-m}.
\end{align*}
\end{cor}
\begin{rem} \label{rem:ajout2}
\begin{enumerate}[(i)] 
\item Note that we have actually proved a more complete result than \eqref{eq:prop:align} and \eqref{bis}:
\begin{equation}
\label{eq:rem:loin}
\mathbb{P}(L_n(1)=k,\, T_{2m}\le n<T_{2m+1})=\widehat{A}_2(m+1,k+1)A_1(m,n-k),
\end{equation}
for $0\le m\le k\wedge (n-k)$
and
\begin{equation}
\label{eq:rem:loin2}
\mathbb{P}(L_n(1)=k,\, T_{2m-1}\le n<T_{2m})=A_2(m,k+1)\widehat{A}_1(m,n-k), 
\end{equation}
for $1\le m\le (k+1)\wedge(n-k)$,
where $\widehat{A}_1$ and $\widehat{A}_2$ are defined by \eqref{eq:defdehatA}.
\item We deduce from \eqref{eq:rem:loin} that
\[
\mathbb{P}(T_{2m}\le n<T_{2m+1})=\sum_{k=m}^{n-m}\widehat{A}_2(m+1,k+1)A_2(m,n-k).
\]
Since the left hand side equals $\mathbb{P}(T_{2m}\le n,\, T_{2m+1}-T_{2m}>n-T_{2m})$, Proposition~\ref{prop:changtpente-sym} and Remark~\ref{rem:2.6} imply
\[
\mathbb{E}[\ind{T_{2m}\le n}\mathcal{P}_2(n-T_{2m})]=\sum_{k=m}^{n-m}\widehat{A}_2(m+1,k+1)A_2(m,n-k).
\]
Recall that $T_{2m}\ge 2m$. Then taking successively $n=2m$, $n=2m+1$ and so on, we are theoreticaly able to determine the law of $T_{2m}$.
\end{enumerate}
\end{rem}

\noindent As it is said in Remark \ref{rem:ajout2}, Proposition \ref{prop:loicm} contains in an hidden way the distribution of $T_{2m}$ and $T_{2m+1}$. However it is actually possible to determine differently the distribution of these two random variables. It is convenient to introduce the notations:
\[
\Delta g(n)=g(n)-g(n+1),\quad n\ge 0,
\]
\[
\varphi*\psi(n)=\sum_{k=0}^n \varphi(k)\psi(n-k)\quad n\ge 0,
\]
\[
\theta:\mathbb{N}\to \mathbb{N}\quad \theta(n)=n+1,
\]
where $g$, $\varphi$, $\psi:\mathbb{N}\to\mathbb{N}$.
\begin{prop}
\label{prop:generali}
Let $\xi_1$,...,$\xi_k$ be $k$ independent $\mathbb{N}$-valued random variables. Denote for any $n \ge 0$, 
\[
f_i(n):=\mathbb{P}(\xi_i\ge n).
\]
We introduce $\mathbb{A}_r^k$ the set of all subsets of $\{1,\ldots,k\}$ containing $r$ elements.
Then 
\[
\mathbb{P}(\xi_1+\ldots+\xi_k\ge n)=h_k(n)
\]
where
\[
h_k=\sum_{r=1}^k\sum_{A\in \mathbb{A}_r^k}\Delta^{r-1}(f^{*A}\circ \theta^{k-r})
\]
and $f^{*A}=f_{i_1}*\ldots * f_{i_r}$ when $A=\{i_1,\ldots,i_r\}$.
\end{prop}
We do not prove Proposition \ref{prop:generali} since it does not play a main role in our study.
\begin{rem}\label{remrem}
\begin{enumerate}
\item If $\xi$ is geometrically distributed with parameter $1-\rho$ (i.e. $\mathbb{P}(\xi=n)=(1-\rho)\rho^n$, $n\ge 1$, $\rho\in]0,1[$) then the function $f$ associated with $\xi$ is $f(n)=\rho^n$, $n\ge 0$.
\item Suppose that $\xi_1=T_{2m}-T_{2m-1}$ (resp. $\xi_2=T_{2m-1}-T_{2m-2}$) where $m\ge 1$, then Remark \ref{rem:2.6} implies that 
\[
\mathbb{P}(\xi_1-1\ge n)=\mathcal{P}_i(n+1),\quad i=1,2,\ n\ge 0,
\]
where $\mathcal{P}_i$ has been defined by \eqref{eq:}
\end{enumerate}
\end{rem}
\begin{defi} Let $\rho\in]0,1[$. A $\mathbb{N}$-valued random variable $\xi_\rho$ is said to be pseudo-Poisson distributed with parameter $\rho>0$ when for all $n\ge 0$:
\[
f_\rho(n)=\mathbb{P}(\xi_\rho\ge n)=\frac{\rho^n}{n!}.
\]
\end{defi}
\noindent It is clear that if $\alpha_{i,k}=1-\frac{\rho_i}{k}$ where $\rho_i\in]0,1[$, then $\mathcal{P}_i(n)=\frac{\rho_i^n}{n!}$. Therefore $\xi_i-1$ (cf item 2. of Remark \ref{remrem}) is pseudo-Poisson with parameter $\rho_i$.\\
It is immediate to prove that:
\[
f_{\rho}*f_{\rho'}=f_{\rho+\rho'}.
\] 
Reasoning by induction on $k$ and using Proposition \ref{prop:generali}, we get the following result.
\begin{prop}\label{prop:pseudopoi}
Suppose that $\xi_1$,...,$\xi_k$ are independent, and $\xi_i$ is pseudo-Poisson with parameter $\rho_i$. Then:
\[
\mathbb{P}(\xi_1+\ldots+\xi_k\ge n)=h_k(n),\quad n\ge 0,
\]
where
\begin{equation}
\label{eq:prop:pois}
h_k(n)=\sum_{r=1}^k\sum_{A\in\mathbb{A}_r^k}\sum_{\ell=0}^{r-1}\binom{r-1}{\ell}\frac{(-1)^\ell}{(n+k+\ell-r)!}\, \left( \sum_{i\in A}\rho_i \right)^{n+k+\ell-r}.
\end{equation} 
In the particular case $\rho_1=\ldots =\rho_k=\rho$,
\[
h_k(n)=k\sum_{t=0}^{k-1}\frac{\rho^{n+t}}{(n+t)!}\binom{k-1}{t}\left( \sum_{\ell=0}^t\binom{t}{\ell}(-1)^\ell(\ell+k-t)^{n-1+t} \right)
\]
\end{prop}
\begin{rem} Suppose that $\alpha_{1,k}=1-\frac{\rho}{k}$, $k\ge 1$. Then 
\[
\mathbb{P}\left(\sum_{i=1}^k(T_{2i}-T_{2i-1}-1)\ge n\right)=h_k(n),\quad n\ge 0,
\]
where $h_k$ is given by \eqref{eq:prop:pois}.
\end{rem}

\subsection{Distribution of the persistent random walk at an independent time}
As shows Proposition \ref{prop:loicm}, the law of $S_n$ is rather complicated. In the study of a Markov chain, it can be interesting to stop it at a random time. For instance, a Markov chain stopped at a geometric time independent from the Markov chain remains a Markov chain.\\ 
Let us consider a geometric random variable $\tau+1$ with parameter $\rho\in]0,1[$ and independent from $(X_n,M_n)$:
\begin{equation}\label{eq:geom}
\mathbb{P}(\tau=k)=\rho^k(1-\rho),\quad k\ge 0.
\end{equation}
In this section we first determine in Theorem \ref{thm:doublelap} below the generating function $\Phi(\lambda,\rho)$ of $S_\tau$:
\begin{equation}
\label{eq:defdedblgen}
\Phi(\lambda,\rho):=\mathbb{E}[\lambda^{S_\tau}]=(1-\rho)\sum_{k\ge 0} \rho^k\mathbb{E}[\lambda^{S_k}],\quad 0\le \rho \le 1.
\end{equation}
This would allow to deduce the generating function of $S_k$ for any $k$ since:
\begin{equation}
\label{eq:defdedblgen2}
\mathbb{E}[\lambda^{S_k}]=\frac{1}{k!}\left.\frac{\partial^k}{\partial \rho^k}\left( \frac{\Phi(\lambda,\rho)}{1-\rho} \right)\right|_{\rho=0}.
\end{equation}
Since we have already calculated the law of $S_k$ we do not go further in this direction.\\
In Section~\ref{sec:GITN}, we will prove that under certain conditions, the persitent random walk $(S_n)$ converges to a Markov process $\left(S(t)\right)_{t\in\mathbb{R}_+}$. The following Theorem~\ref{thm:doublelap} will be used to calculate the Laplace transform of $S(\xi)$, where $\xi$ is an exponential random variable independent of $\left(S(t)\right)_{t\in\mathbb{R}_+}$. Theoretically, the following theorem permits to deduce the law of $S_\tau$ but it is in practice impossible to determine it explicitely. However, using the law of $L_n(1)$ for any $n$, given in Proposition~\ref{prop:loicm}, we present in Proposition~\ref{prop:ntau} below the distribution of $L_\tau(1)$. Recall that from \eqref{eq:prop:loicm1}, $S_\tau=1+2L_\tau(1)-\tau$. Since $\tau$ is a random time, we cannot deduce from this indentity the distribution of $S_\tau$.
\begin{thm}
\label{thm:doublelap}
Let $0<\rho<\lambda<1$. Then the generating function of $S_\tau$, where $(S_n)_{n\geq 0}$ and $\tau$ are independent, is equal to
\begin{equation}
\label{eq:thm:double}
\mathbb{E}[\lambda^{S_\tau}]=\frac{(\rho-1)\left\{\lambda\rho\left(\widehat{\mathcal{P}}_1\left(\frac{\rho}{\lambda}\right)+\widehat{\mathcal{P}}_2(\lambda\rho)\right)+(\lambda\rho-1)\widehat{\mathcal{P}}_1\left(\frac{\rho}{\lambda}\right)\widehat{\mathcal{P}}_2(\lambda\rho)\right\}}{\rho(\lambda\rho-1)\widehat{\mathcal{P}}_2(\lambda\rho)+\lambda\rho(\rho-\lambda)\widehat{\mathcal{P}}_1\left(\frac{\rho}{\lambda}\right)+(\lambda\rho-1)(\rho-\lambda)\widehat{\mathcal{P}}_1\left(\frac{\rho}{\lambda}\right)\widehat{\mathcal{P}}_2(\lambda\rho)}
\end{equation}
where $\widehat{\mathcal{P}}_i$ is defined for $ 0<x< 1$ by 
\begin{equation}
\label{eq:def:pchap}
\widehat{\mathcal{P}}_i(x)=\sum_{k\ge 1} \mathcal{P}_i(k)x^k,\quad  i=1,2.
\end{equation}
\end{thm}
%We focus our attention to an example.
\begin{rem} If $\alpha_{2,k}=\alpha_2$ for any $k\ge 1$ (recall that in that case $S_n$ is the persistent random walk associated with the simple infinite comb), the function $\widehat{\mathcal{P}}_2$ satisfies
\[
\widehat{\mathcal{P}}_2(x)=\sum_{k\ge 1}(1-\alpha_2)^{k-1}x^k=\frac{x}{1-(1-\alpha_2)x}.
\]
Therefore \eqref{eq:thm:double} becomes
\[
\mathbb{E}[\lambda^{S_\tau}]=\frac{\lambda(\rho-1)\left(1-\alpha_2\widehat{\mathcal{P}}_1\left(\frac{\rho}{\lambda}\right)\right)}{\lambda\rho-1+\alpha_2\lambda(\rho-\lambda)\widehat{\mathcal{P}}_1\left(\frac{\rho}{\lambda}\right)}.
\]
Moreover, if $\alpha_{1,k}=1-\alpha_1/k$, then
\[
\widehat{\mathcal{P}}_1(x)=\sum_{k\ge 1}\frac{\alpha_1^{k-1}x^k}{(k-1)!}=xe^{\alpha_1 x},
\]
and
\[
\mathbb{E}[\lambda^{S_\tau}]=\frac{(\rho-1)\left(\lambda-\alpha_2\rho e^{\alpha_1\rho/\lambda}\right)}{\lambda\rho-1+\alpha_2\rho(\rho-\lambda)e^{\alpha_1\rho/\lambda}}.
\]
\end{rem}
\noindent We begin with a preliminary result (Lemma \ref{lem:reltop}). The proof of Theorem \ref{thm:doublelap} will be given later on.  For $i\in \{1,2\}$ and $0<x<1$, let us define the generating function 
\[\mathcal{G}^{(i)}(x):=\sum_{k\ge 1}\mathcal{P}_i(k)\alpha_{i,k} x^k.\]
\begin{lem}\label{lem:reltop}
\begin{enumerate}[(i)]
\item For $i=1,2$ and $0<x<1$, the generating function $L^{(i)}$ satisfies
\begin{equation}\label{eq:reltop}
\mathcal{G}^{(i)}(x)=1+\left(\frac{x-1}{x}\right)\widehat{\mathcal{P}}_i(x),
\end{equation}
where $\widehat{\mathcal{P}}_i(x)$ has been defined by \eqref{eq:def:pchap}.
\item Moreover for $m\ge 1$, 
\begin{equation}
\label{eq:lem111}
\sum_{b\ge m}A_i(m,b)x^b=\left(\mathcal{G}^{(i)}(x)\right)^m,\quad \sum_{b\ge m}\widehat{A}_i(m,b)x^b=\left(\mathcal{G}^{(i)}(x)\right)^{m-1}\widehat{\mathcal{P}}_i(x),
\end{equation}
where $A_i$ (resp. $\widehat{A}_i$) is defined by \eqref{eq:defdeA} (resp. \eqref{eq:defdehatA}).
\end{enumerate}
\end{lem}
{\noindent {\sc Proof of Lemma \ref{lem:reltop}.}\ }
\begin{enumerate}[(i)]
\item Let $0<x<1$. We have
\begin{align*}
\mathcal{G}^{(i)}(x)&=-\sum_{k\ge 1}\mathcal{P}_i(k)(1-\alpha_{i,k})x^k+\sum_{k\ge 1}\mathcal{P}_i(k)x^k\\
&=-\frac{1}{x}\, \sum_{k\ge 1}\mathcal{P}_i(k+1)x^{k+1}+\widehat{\mathcal{P}}_i(x)=-\frac{1}{x}(\widehat{\mathcal{P}}_i(x)-x)+\widehat{\mathcal{P}}_i(x).
\end{align*}
\item For $m\ge 1$,
\begin{align*}
\sum_{b\ge m}A_i(m,b)x^b&= \sum_{b\ge m,\ u\in\mathcal{N}(m,b)}\mathcal{P}_i(u_1)\ldots \mathcal{P}_i(u_m)\ \alpha_{i,u_1}\times\ldots\alpha_{i,u_m}x^{u_1+\ldots+u_m}\\
&=\sum_{u\in(\mathbb{N}^*)^m}\Big( \mathcal{P}_i(u_1)\ \alpha_{i,u_1} x^{u_1} \Big)\ldots\Big( \mathcal{P}_i(u_m)\ \alpha_{i,u_m} x^{u_m} \Big)\\
&=\left(\mathcal{G}^{(i)}(x)\right)^m.
\end{align*}
\noindent The proof of the second equality in \eqref{eq:lem111} is similar to the first one.\QED
\end{enumerate}
{\noindent {\sc Proof of Theorem~\ref{thm:doublelap}.}\ }
Let $0<\rho<\lambda<1$. Using \eqref{eq:prop:loicm1} together with the independence between $S_n$ and $\tau$ yield
\begin{align}
\label{eq:thm:eq1}
\mathbb{E}[\lambda^{S_\tau}]&=(1-\rho)\sum_{n\ge 0}\mathbb{E}[\lambda^{S_n}]\rho^n=\lambda(1-\rho)\sum_{n\ge 0}\mathbb{E}[\lambda^{2L_n(1)}]\left(\frac{\rho}{\lambda}\right)^n\nonumber\\
&=\lambda(1-\rho)\sum_{k\ge 0}\lambda^{2k}\sum_{n\ge k}\eta_n(k)\left(\frac{\rho}{\lambda}\right)^n.
\end{align}
See Proposition \ref{prop:loicm} for the definition of $\eta_k(n)$.
Using the decomposition \eqref{eq:proof} and equality \eqref{eq:etap11} lead to the following decomposition
\begin{equation}
\label{eq:dec:thm}
\mathbb{E}[\lambda^{S_\tau}]=\lambda(1-\rho)(\mathcal{E}_1+\mathcal{E}_2),
\end{equation}
where $\mathcal{E}_i$ corresponds to the part related to $\eta_n^{(i)}$ (cf \eqref{eq:prop:align} and \eqref{bis})  i.e. :
\begin{align*}
\mathcal{E}_1&=\sum_{k\ge 0}\lambda^{2k}\sum_{n\ge k}\sum_{m=0}^{k\wedge (n-k)}\eta_n^{(1)}(k,m)\left(\frac{\rho}{\lambda}\right)^n\\
&=\sum_{k\ge 0}(\lambda\rho)^{k}\sum_{m=0}^{k}\widehat{A}_2(m+1,k+1)\sum_{n\ge m+k}A_1(m,n-k)\left(\frac{\rho}{\lambda}\right)^{n-k}.
\end{align*}
By \eqref{eq:lem111}, we get
\begin{align*}
\mathcal{E}_1&=\sum_{k\ge 0}(\lambda\rho)^{k}\sum_{m=0}^{k}\widehat{A}_2(m+1,k+1)\left(\mathcal{G}^{(1)}\Big(\frac{\rho}{\lambda}\Big)\right)^{m}\\
&=\frac{1}{\lambda\rho}\sum_{m\ge 0}\left(\mathcal{G}^{(1)}\Big(\frac{\rho}{\lambda}\Big)\right)^{m}\sum_{k\ge m}\widehat{A}_2(m+1,k+1)(\lambda\rho)^{k+1}\\
&=\frac{1}{\lambda\rho}\sum_{m\ge 0}\left(\mathcal{G}^{(1)}\Big(\frac{\rho}{\lambda}\Big)\mathcal{G}^{(2)}(\lambda\rho)\right)^{m}\widehat{\mathcal{P}}_2(\lambda\rho)=\frac{\widehat{\mathcal{P}}_2(\lambda\rho)}{\lambda\rho\Big(1-\mathcal{G}^{(1)}\Big(\frac{\rho}{\lambda}\Big)\mathcal{G}^{(2)}(\lambda\rho)\Big)}.
\end{align*}
In a similar way, we compute $\mathcal{E}_2$:
\begin{align*}
\mathcal{E}_2&=\sum_{k\ge 0}\lambda^{2k}\sum_{n\ge k}\sum_{m=1}^{(k+1)\wedge (n-k)}\eta_n^{(2)}(k,m)\left(\frac{\rho}{\lambda}\right)^n\\
&=\sum_{k\ge 0}(\lambda\rho)^{k}\sum_{m=1}^{k+1}A_2(m,k+1)\sum_{n\ge m+k}\widehat{A}_1(m,n-k)\left(\frac{\rho}{\lambda}\right)^{n-k}\\
&=\frac{1}{\lambda\rho}\ \widehat{\mathcal{P}}_{1}\left(\frac{\rho}{\lambda}\right)\sum_{m\ge 1}\left(\mathcal{G}^{(1)}\Big(\frac{\rho}{\lambda}\Big)\right)^{m-1}\sum_{k\ge m-1}A_2(m,k+1)(\lambda\rho)^{k+1}\\
&=\frac{\widehat{\mathcal{P}}_{1}\left(\frac{\rho}{\lambda}\right)\mathcal{G}^{(2)}(\lambda\rho)}{\lambda\rho\Big(1-\mathcal{G}^{(1)}\Big(\frac{\rho}{\lambda}\Big)\mathcal{G}^{(2)}(\lambda\rho)\Big)}.
\end{align*}
Now \eqref{eq:dec:thm} yields
\begin{align*}
\mathbb{E}[\lambda^{S_\tau}]&=\frac{\lambda(1-\rho)}{\lambda\rho}\frac{\widehat{\mathcal{P}}_2(\lambda\rho)+\widehat{\mathcal{P}}_{1}\left(\frac{\rho}{\lambda}\right)\mathcal{G}^{(2)}(\lambda\rho)}{1-\mathcal{G}^{(1)}\left(\frac{\rho}{\lambda}\right)\mathcal{G}^{(2)}(\lambda\rho)}
\end{align*}
which, combined with \eqref{eq:reltop},  implies \eqref{eq:thm:double}. \QED
 %where $\tau$ satisfies \eqref{eq:geom}.
\begin{prop}\label{prop:ntau}
Let $k\ge 0$. The random variable $L_\tau(1)$ satisfies 
\begin{align}
\label{eq:ntau1}
\mathbb{P}(L_\tau(1)=k)&=(1-\rho)\rho^k\Big\{g_1(\rho)\sum_{m=1}^{k+1} A_2(m,k+1)f_1(\rho)^{m-1}\nonumber\\
&+\sum_{m=0}^{k}\sum_{\ell=1}^{k-m+1}A_2(m,k+1-\ell)\mathcal{P}_2(\ell)f_1(\rho)^m\Big\}
\end{align}
with
\begin{equation}
\label{eq:ntau2}
f_i(\rho)=\sum_{k\ge 1}\mathcal{P}_i(k)\alpha_{i,k}\rho^k\quad\mbox{and}\quad
g_i(\rho)=\sum_{k\ge 1}\mathcal{P}_i(k)\rho^k,\quad i=1,2.
\end{equation}
Moreover 
\begin{equation}
\label{eq:gifi}
f_i(\rho)=\left(1-\frac{1}{\rho}\right)g_i(\rho)+1.
\end{equation}
\end{prop}
\pff
Let us first recall (cf Proposition \ref{prop:loicm}) that 
\[
\eta_n(k):=\mathbb{P}(L_n(1)=k)=\eta_n^{(1)}(k)+\eta_n^{(2)}(k),
\]
where $\eta_n^{(1)}$ resp. $\eta_n^{(2)}$ is defined by \eqref{eq:prop:align} resp. \eqref{bis}. In a similar way, we decompose the following probability
\begin{equation}\label{sandou}
\eta(k):=\mathbb{P}(L_\tau(1)=k)=\eta^{(1)}(k)+\eta^{(2)}(k)
\end{equation}
where 
\begin{align*}
\eta^{(i)}(k)=(1-\rho)\sum_{n\ge 0}\rho^n\eta_n^{(i)}(k),\quad i=1,2. 
\end{align*}
We shall only present the details of calculation for $\eta^{(1)}(k)$ ($\eta^{(2)}(k)$ can be determined similarly). By definition
\begin{equation}\label{eq:devlpterm1}
\eta^{(1)}(k)=(1-\rho)\sum A_2(m,k+1)A_1(m-1,n-k-\ell)\mathcal{P}_1(\ell)\rho^n,
\end{equation}
the sum is taken over all combinations of indexes $n$, $m$,  and $\ell$ satisfying
\[
n\ge 0,\quad m\ge 1, \quad m\le k+1,\quad m\le n-k,\quad \ell\ge 1,\quad \ell\le n-k-m+1.
\]
Let us first fix the indexes $m$ and $\ell$ with
\begin{equation}
\label{eq:domaindef}
1\le m,\quad m\le k+1,\quad \ell\ge 1.
\end{equation}
Then we compute the sum with respect to $n$. We therefore introduce
\[
\psi_{m,\ell}(n):=\sum_{n\ge \ell+k+m-1}A_1(m-1,n-k-\ell)\rho^n.
\]
By the change of variable $i=n-k-\ell-m+1$, we get 
\begin{align*}
\psi_{m,\ell}(n)&=\sum_{i\ge 0}A_1(m-1, m-1+i)\rho^{k+\ell}\rho^{m-1+i}\\
&=\rho^{k+\ell}\sum_{i\ge 0}\sum_{u_1,\ldots, u_{m-1}}\mathcal{P}_1(u_1)\times\ldots\times \mathcal{P}_1(u_{m-1})\alpha_{1,u_1}\times\ldots\times\alpha_{1,u_{m-1}}\\
&\times \rho^{u_1+\ldots+u_{m-1}}
\ind{u_1+\ldots+u_{m-1}=m-1+i}\\
&=\rho^{k+\ell}\sum_{u_1,\ldots, u_{m-1}}\mathcal{P}_1(u_1)\times\ldots\times \mathcal{P}_1(u_{m-1})\alpha_{1,u_1}\times\ldots\times\alpha_{1,u_{m-1}}\rho^{u_1+\ldots+u_{m-1}}\\
&=\rho^{k+\ell}(f_1(\rho))^{m-1},
\end{align*}
where $f_1$ is defined by \eqref{eq:ntau2}. Let us just note that, in the particular case $m=1$, we get $A_1(0,n-k-\ell)=\ind{n-k-\ell=0}$ and $\psi_{1,\ell}(k)=\rho^{k+\ell}$. Using \eqref{eq:devlpterm1} we obtain the following sum over all indexes $m$ and $\ell$ satisfying \eqref{eq:domaindef}:
\begin{align*}
\eta^{(1)}(k)=(1-\rho)\sum A_2(m,k+1)\mathcal{P}_1(\ell)\rho^{k+\ell}\left(f_1(\rho)\right)^{m-1},
\end{align*}
when $m,\ell$ verify \eqref{eq:domaindef}. Then
\begin{align}\label{eq:fin}
\eta^{(1)}(k)=(1-\rho)\rho^k g_1(\rho)\left( \sum_{m=1}^{k+1} A_2(m,k+1) f_1(\rho)^{m-1} \right)
\end{align}
where $g_1$ is defined by \eqref{eq:ntau2}. \\
It can be proved 
\begin{equation}
\label{eq:fin2}
\eta^{(2)}(k)=(1-\rho)\rho^k \sum_{m=0}^{k}\sum_{\ell= 1}^{ k-m+1}A_2(m,k+1-\ell)\mathcal{P}_2(\ell)f_1(\rho)^m.
\end{equation}
Obviously \eqref{sandou}, \eqref{eq:fin} and \eqref{eq:fin2} imply \eqref{eq:ntau1}.
Let us finally prove \eqref{eq:gifi}:
\begin{align*}
f_i(\rho)&=\sum_{k\ge 1}\mathcal{P}_i(k)\alpha_{i,k}\rho^k=\sum_{k\ge 1}\mathcal{P}_i(k)\left(1-(1-\alpha_{i,k})\right) \rho^k=g_i(\rho)-\sum_{k\ge 1}\mathcal{P}_i(k+1)\rho^k\\
&=g_i(\rho)-\frac{1}{\rho}\sum_{k\ge 2}\mathcal{P}_i(k)\rho^k=g_i(\rho)-\frac{1}{\rho}(g_i(\rho)-\rho)=\Big( 1-\frac{1}{\rho} \Big)g_i(\rho)+1.
\quad \blacksquare \end{align*}

\subsection{Large time behavior}\label{sec:largetime}
The law of $S_n$ has been given explicitely in Proposition~\ref{prop:loicm} but it is very complicated. This leads us to investigate the asymptotic behaviour of $S_n$ as $n\to\infty$. 
\begin{prop}
\label{prop:asympt}
Assume that $\Theta_i<\infty$, $i=1,2$, where $\Theta_i$ is defined by \eqref{eq:cond-gene}.
\begin{enumerate}[(i)] 
\item  The ratio $\displaystyle \frac{S_n}{n}$ converges a.s. and in $L^1$ to $\displaystyle\frac{\Theta_2-\Theta_1}{\Theta_1+\Theta_2}$ as $n\to\infty$.
\item Moreover, if $\sum_{k\ge 1}k\mathcal{P}_i(k)<\infty$ for $i=1,2$, then the Central Limit Theorem holds:
\begin{equation}
\label{eq:TCL}
\frac{1}{\sqrt{n} \Upsilon}\left(S_n- n\frac{\Theta_2-\Theta_1}{\Theta_1+\Theta_2} \right)
\end{equation}
converges in distribution to a standard Gaussian random variable as $n\to\infty$ and the constant $\Upsilon$ is defined by
\begin{equation}
\label{eq:upsilon}
\Upsilon=\frac{4}{\Theta_1+\Theta_2}\ \mathbb{E}\left[\left( T_1-\frac{\Theta_2 T_2}{\Theta_1+\Theta_2} \right)^2\right]
\end{equation}
where the stopping times $T_1$ and $T_2$ are defined by \eqref{eq:changt} and $X_0=M_0=1$.
\end{enumerate}
\end{prop}
\begin{rem}
\begin{enumerate}
\item Let us first note that, under the condition presented in (ii) we can also prove the existence of a constant $C\in\mathbb{R}$ such that 
\begin{equation}
\label{eq2:prop:asympt}
\lim_{n\to\infty}\left \{ \mathbb{E}(S_n)-n\ \frac{\Theta_2-\Theta_1}{\Theta_1+\Theta_2} \right\}=C.
\end{equation}
\item In the particular case $\Theta_1=\Theta_2<\infty$, Proposition~\ref{prop:asympt} implies that $\lim_{n\to\infty}\frac{\mathbb{E}(S_n)}{n}=0$. If moreover $\sum_{k\ge 0} k \mathcal{P}_i(k)<\infty$, we have a more precise result which says that $\frac{1}{\sqrt{n}}S_n$ converges in distribution to a Gaussian random variable.
\item Under the conditions $\Theta_i<\infty$ and $\sum_{k\ge 1}k\mathcal{P}_i(k)<\infty$, we observe therefore that the rates of convergence for the first and the second order limit theorems are similar to the rates in the setting of the classical Bernoulli random walk. The persistency does not change the long time behaviour. 
\item The assumption $\sum_{k\ge 1}k\mathcal{P}_i(k)<\infty$ is quite strong and force a relatively strong mixing in the sequence $(X_n)$. Open and interesting questions occur when this assumption is not satisfied. In terms of VLMC, it corresponds to the case when the expectation of the length of $\lpref(U_n)$ is infinite. 
\end{enumerate}
\end{rem}
{\noindent {\sc Proof of Proposition \ref{prop:asympt}}\ }
\begin{enumerate}[(i)]
\item Proposition~\ref{prop:stat-gen} ensures that, under the condition $\Theta_i<\infty$, for $i\in \{1,2\}$, the process $(X_n,M_n)_{n\ge 0}$ is an ergodic Markov chain with invariant probability $\nu$. The ergodic theorem, Corollary~\ref{cor:stat-gen} and \eqref{eq:cond-gene} imply the following almost sure convergence result:
\begin{equation}\label{eq:proof:asym}
\lim_{n\to\infty}\frac{L_n(1)}{n}=\nu(1,\mathbb{N})=\frac{\Theta_2}{\Theta_1+\Theta_2}\quad a.s.,
\end{equation}
where $L_n(1)$ is defined by \eqref{eq:prop:loicm}. Since $L_n(1)/n$ is a bounded random variable, the almost sure convergence implies the moment convergence. Therefore, by \eqref{eq:prop:loicm1} and \eqref{eq:proof:asym}, we obtain
\[
\lim_{n\to\infty}\frac{\mathbb{E}(S_n)}{n}=\lim_{n\to\infty}\frac{1+2\mathbb{E}(L_n(1))}{n}-1=\frac{2\Theta_2}{\Theta_1+\Theta_2}-1=\frac{\Theta_2-\Theta_1}{\Theta_1+\Theta_2}.
\]
\item Let us consider the Markov chain $(X_n,M_n)_{n\ge 0}$ starting at $(1,1)$ and denote $Q$ the associated transition probability and $\nu$ the invariant measure. We define 
\begin{equation}\label{eq:def:sig}
\sigma=\inf\Big\{ n\ge 1:\ (X_n,M_n)=(1,1) \Big\}.
\end{equation}
Since the Markov chain is reccurent irreducible and positive, the stopping time $\sigma$ is almost surely finite. Moreover if $\mathbb{E}[\sigma^2]<\infty$, Theorem~17.2.2 in \cite{Meyn-Tweedie} implies that \eqref{eq:TCL} holds with the constant
\begin{equation}\label{eq:ups}
\Upsilon:=\nu(1,1)\ \mathbb{E}\left[ \left( \sum_{k=1}^\sigma \Big(X_k-\frac{\Theta_2-\Theta_1}{\Theta_1+\Theta_2}\Big)\right)^2 \right].
\end{equation}
%
%
% implies that 
%\begin{equation}\label{cauchy}
%\sum_{n\ge 0}d_{TV}\left(\delta_{(1,1)}Q^n,\nu\right)<\infty,
%\end{equation}
%where $\delta_{(1,1)}$ is the Dirac measure at $(1,1)$ and the total variation distance is defined by
%\[
%d_{TV}(\mu,\nu)=\sup_{A\subset \{-1,1\}\times \mathbb{N}^*}|\mu(A)-\nu(A)|.
%\]
%Let us assume now that $\mathbb{E}[\sigma^2]<\infty$ and define
%\begin{align*}
%u_n&=\mathbb{E}[S_n]-n\frac{\Theta_2-\Theta_1}{\Theta_1+\Theta_2}=1+2\mathbb{E}[L_n(1)]-n-n\frac{\Theta_2-\Theta_1}{\Theta_1+\Theta_2}\\
%&=1+2\left( \mathbb{E}[L_n(1)]-\frac{n\Theta_2}{\Theta_1+\Theta_2}\right).
%\end{align*}
%Therefore
%\begin{align*}
%|u_{n+1}-u_n|&=2\left|\mathbb{P}(X_{n+1}=1)-\frac{\Theta_2}{\Theta_1+\Theta_2}\right|\\
%&=2|\mathbb{P}(X_{n+1}=1)-\nu(1,\mathbb{N}^*)|\le d_{TV}\left(\delta_{(1,1)}Q^{n+1},\nu\right).
%\end{align*}
%Finally \eqref{cauchy} implies that $(u_n)_{n\geq 0}$ is a Cauchy sequence in $\mathbb{R}$. Consequently $(u_n)$ converges towards a limit $C$, as $n\to\infty$. \\
According to Definition \eqref{eq:changt} of the stopping times $(T_n)$, one has $\sigma=T_2$ and consequently
\[
\sum_{k=1}^\sigma X_k=\sum_{k=1}^{T_2-1}X_k+\sum_{k=T_1}^{T_2-1}X_k+X_{T-2}=T_1-1-(T_2-T_1)+1=2T_1-T_2.
\]
From \eqref{inv-mes-double} and \eqref{eq:ups}, we deduce \eqref{eq:upsilon}. It remains to prove that $\sigma$ is square integrable. Since $\sigma=T_1+(T_2-T_1)$ and $T_2-T_1\ge 0$, $\mathbb{E}(\sigma^2)<\infty$ if and only if $\mathbb{E}[T_1^2]<\infty$ and $\mathbb{E}[(T_2-T_1)^2]<\infty$. Using Proposition~\ref{prop:changtpente-sym} we have: 
\begin{align*}
\mathbb{E}[T_1^2]&=\sum_{n\ge 1}n^2\mathcal{P}_2(n)\ \alpha_{2,n}=-\lim_{N\to\infty}\sum_{n= 1}^N n^2\mathcal{P}_2(n)\ \left((1-\alpha_{2,n})-1\right)\\
&=-\lim_{N\to\infty}\left(\sum_{n=1}^N n^2\mathcal{P}_2(n+1)-\sum_{n=1}^N n^2\mathcal{P}_2(n)\right)\\
&\le\mathcal{P}_2(1)+\lim_{N\to\infty}\sum_{n=2}^N(n^2-(n-1)^2)\mathcal{P}_2(n)\\
&\le 1+ \lim_{N\to\infty}\sum_{n=2}^N(2n-1)\mathcal{P}_2(n)\le 1+2\sum_{n\ge 1}n\mathcal{P}_2(n)<\infty.
\end{align*}
Using \eqref{eq:dist++} and similar arguments, we obtain that $\mathbb{E}[(T_2-T_1)^2]<\infty$. \QED
\end{enumerate}

\section{From persistent random walk to generalized integrated telegraph noise (GITN).}
\label{sec:GITN}
Let $(X_n,M_n)_{n\ge 0}$ be a $\{-1,1\}\times\overline{\mathbb{N}}^*$-valued Markov chain satisfying \eqref{eq:transition} and \eqref{eq:extension} and let $(S_n)_{n\ge 0}$ be the associated persistent random walk, see \eqref{eq:def-persist-part}. We assume in this section that the transition probabilities $(\alpha_{i,n})$ depend on a small parameter $\eps>0$ and $\eps$ appears also both in a time scale and a space scale of the persistent random walk. We prove that there exists a normalization expressed in terms of $\eps$ so that $(X_n,M_n,S_n)$ converges in distribution as $\eps\to 0$. This limit is a time continuous process. Such a procedure has been already performed in \cite{Herrmann-Vallois} when the increments are a Markov chain.\\
More precisely we suppose that the transition probabilities satisfy
\begin{equation}
\label{eq:assum-asym}
\alpha_{i,n}=f_i(n\eps)\eps+\tilde{\alpha}_{i,n,\eps}\eps, \quad n\ge 1,\ i=1,2
\end{equation}
where $f_1$ and $f_2$ are positive right-continuous functions with left limits satisfying
\begin{equation}
\label{eq:int-fini} \int_0^\infty f_i(u)du=\infty, \ i=1,2
\end{equation}
and $\tilde{\alpha}_{i,n,\eps}\in\mathbb{R}$ with $\lim_{\eps\to 0}\sup_{i,n}|\tilde{\alpha}_{i,n,\eps}|=0$.
It is clear that for any $i,n$ fixed, $\lim_{\eps\to 0}\alpha_{i,n}=0$. Therefore $X_k$ changes from $-1$ to $1$ (for instance) with a small probability. The trend of $(X_k)$ is to stay at the same level. 
\\  Let us now introduce the scaling procedure. For any $\eps >0$ and for any $t\in \eps\mathbb{N}$, we
define the processes
\begin{equation}\label{eq:defdeseps}
S^\eps(t)=\eps S_{\frac t{\eps}},\quad M^\eps(t)=\eps M_{\frac t{\eps}}\quad\mbox{and}\quad X^\eps(t)=X_{\frac t{\eps}}.
\end{equation}
%where $(S_n)_{n\geq 0}$ is defined by \eqref{eq:def-persist-part}. 
Note that $(S_n)$ depends on $\eps$, since the two families of coefficients $(\alpha_{1,n})$ and $(\alpha_{2,n})$ depend on $\eps$. For the sake of simplicity, we do not mention the dependency with respect to $\eps$.
We extend the definition of the process $(S^\eps(t),\, t\in \eps\mathbb{N})$ to $t\in \mathbb{R}_+$ by linear interpolation and we the definition of the processes $(X^\eps(t),\, t\in \eps\mathbb{N})$ and $(M^\eps(t),\, t\in \eps\mathbb{N})$ into piecewise constant right continuous with left limits functions. In order to describe the asymptotic behavior of $(S^\eps(t),\,t\ge0)$ as $\eps\to 0$, it suffices to study the asymptotic properties of the times of trend changes. Indeed $t\to S^\eps(t)$ admits a $1$ slope till the stopping time $\eps T_1$, with $T_1$ defined by \eqref{eq:changt}.
After that instant, the paths admits a $-1$ slope till $\eps T_2$ and so on... The increments change periodically from $-1$ to $1$ and vice versa. \\
As $\eps\to 0$, we shall prove that the limit process $(S^0(t),\ t\ge 0)$ is still piecewise linear.
More precisely it starts at $t=0$ with a slope equal to $1$. At a random time time $e_1$ the slope changes and becomes equal to $-1$,
at random time $e_1+e_2$ we observe a new change of slope and so on... We are therefore particularly interested in the description of the distribution of $(e_n)_{n\ge 1}$.
\begin{thm}\label{thm:convergenceGITN}
1. Let us consider a sequence $(e_n)_{n\ge 1}$ of independent random variables such that for $n\geq 1$,
\begin{equation}\label{eq:conv-fin}
\mathbb{P}(e_{2n-1}>t)=\exp\left(-\int_0^t f_2(u)du\right),\quad \mathbb{P}(e_{2n}>t)=\exp\left(-\int_0^tf_1(u)du\right),
\end{equation}
where $f_1$ and $f_2$ have been introduced in \eqref{eq:assum-asym}.
Let
\[
N^0(t):=\sum_{n\ge 1}\ind{e_1+\ldots+e_n\le t},\quad \mbox{for any} \ t\ge 0
\]
be the \emph{counting process}, 
\[
m(t):=t-\sup\{e_1+\ldots+e_k:\ e_1+\ldots+e_k\le t\}=t-T_{N^0(t)}
\]
the associate \emph{age process} (spent life)
and finally 
\begin{equation}\label{def:gitn}
S^0(t)=\int_0^t(-1)^{N^0(s)}\,ds,\quad t\ge 0.
\end{equation}
the so-called \emph{Generalized Integrated Telegraph Noise (GITN)}.\\
2. Let $(X_n,M_n)_{n\ge 0}$ be a $\{-1,1\}\times\mathbb{N}^*$-valued Markov chain whose probability transition satisfies \eqref{eq:transition} and is $\eps$-dependent in the sense of \eqref{eq:assum-asym}. We assume $X_0=M_0=1$.
\begin{enumerate}
\item[(i)] For all $n\ge 1$, the sequence of times between two consecutive slope changes $(\eps T_1,\eps (T_2-T_1),\ldots, \eps (T_n-T_{n-1}))$ converges in distribution towards $(e_1,\ldots,e_n)$ as $\eps\to 0$, where the  sequence $(T_k)_{k\ge 0}$ is defined by \eqref{eq:changt}.
\item[(ii)] The following convergence in distribution in Skorohod's topology holds
\begin{equation}
\label{eq:skoro}
\left( S^\eps(t), X^\eps(t), M^\eps(t), t\ge 0  \right)\underset{\eps\to 0}{\longrightarrow} \left( S^0(t), (-1)^{N^0(t)}, m(t), t\ge 0  \right),
\end{equation}
where $S^\eps(t)$, $M^\eps(t)$ and $X^\eps(t)$ are defined by \eqref{eq:defdeseps}.\\
Moreover $\left( S^0(t), (-1)^{N^0(t)}, m(t), t\ge 0  \right)$ and $\left((-1)^{N^0(t)}, m(t), t\ge 0  \right)$ are
Markov processes.
\end{enumerate}
\end{thm}
\begin{rem}\begin{enumerate}\item[(i)]
In the case $X_0=-1$, the family of processes $(S^\eps(t))_{t\ge0}$ converges
in distribution to $(S^0(t))_{t\ge 0}$ as $\eps$ goes to zero, where for any $t\ge 0$,
\[
S^0(t)=-\int_0^t(-1)^{\tilde{N}^0(s)} ds,\quad \mbox{and}\quad \tilde{N}^0(t)=\sum_{n\ge 1}\ind{e_1+\ldots+e_{n+1}\le t}.
\]
In the particular case where the functions $f_1$ and $f_2$ are constant, it has been proved in \cite{Herrmann-Vallois} that a particular solution of the telegraph equation can be represented in terms of $S^0(t)$. That explains that $(S^0(t))$ defined by \eqref{def:gitn} is called the Generalized Integrated Telegraph Noise (GITN).
\item[(ii)]
In the classical integrated telegraph noise
\cite{Herrmann-Vallois}, the random variables $(e_n,\, n\ge 0)$ are exponentially distributed, therefore
$(S^0(t),N^0(t))$ is Markovian. For the generalized situation, this property is not true anymore,
we need to consider some additional information. This information is given by $D_-$ the left derivate of the GITN which is directly related to the age process 
\[
m(t)=t-\sup\{s\ge 0:\ D_-S^0(s)\neq D_-S^0(t)\}.
\]
\item[(iii-a)] Davis wrote in \cite{Davis} that "almost all the continuous-time stochastic process models of applied probability consist of some combination of the following: diffusion, deterministic motion and random jumps". According to Theorem \ref{thm:convergenceGITN}, between two consecutive random jumps the GITN moves in a deterministic way and therefore  belongs to the family of the so-called \emph{Piecewiese Deterministic Markov Processes}, see for instance \cite{Davis,Davisbook,Cocozza}.
\item[(iii-b)] The possible values of $X^0(t)$ are $\{-1,1\}$. It is possible to deal with the case where $X^0(t)\in\{a_1,\ldots, a_K \}$. In that case $X^0(t)$ is a Markov chain indexed by $\mathbb{R}_+$ and $\{a_1,\ldots,a_K\}$-valued. This situation has been already treated in \cite{Herrmann-Vallois}, when the functions $(f_i)_{1\le i\le K}$ are constant.
\item[(iii-c)] $(S^0(t);\ t\ge 0)$ is a semi-Markov process, see \cite{cinlar, janssen2006}. In \cite{janssen2006} (Theorem 3.3 in Chapter 4) it has been proved that $(X^\eps(t);\ t\ge 0)$ converges to the semi-Markov process $(X^0(t),\ t\ge 0)$. This result is weaker than ours since we have considered the convergence of $(S^\eps(t),M^\eps(t),X^\eps(t))_{t\ge 0}$.
\end{enumerate}\end{rem}
%\begin{rem} If the function $f$ satisfies $\int_0^\infty f(u)\,du<\infty$, then $\mathbb{P}(e_{2n}=\infty)>0$ for all $n\ge 0$ and therefore $\sup_{t\ge 0}N^0(t)$ is finite a.s.
%\end{rem}
\pff\paragraph{Step 1 --- Convergence of the jump times.} Let us define $\mathcal{R}^\eps_n:=(\eps T_1, \eps T_2,\ldots,\eps T_n)$ for $n\ge 1$. According to Proposition~\ref{prop:changtpente-sym}, $(T_{n}-T_{n-1})_{n\ge 1}$ is a sequence of independent random variables. In order to prove the convergence in distribution of $\mathcal{R}^\eps_n$ as $\eps$ tends to $0$, it suffices to analyze the behaviour of $\eps(T_{n}-T_{n-1})$ where $n\ge 0$ is given. Recall that $T_0=0$. Remark~\ref{rem:2.6} and \eqref{eq:dist++} yield:
\begin{align*}
\mathbb{P}\left(\eps (T_{2n+1}- T_{2n})>t\right)&=\mathbb{P}\left(T_{2n+1}- T_{2n}>\frac{t}{\eps}\right)=\mathbb{P}\left(T_{2n+1}- T_{2n}> \left\lfloor \frac{t}{\eps}\right\rfloor\right)\\
&=(1-\alpha_{2,1})\times\ldots\times (1-\alpha_{2,\lfloor \frac{t}{\eps}\rfloor}),
\end{align*}
where $\lfloor a \rfloor$ stands for the integer part of $a$.
Defining
\[
\delta_\eps(t):=\log\Big\{\mathbb{P}\Big(\eps( T_{2n+1}- T_{2n})>t\Big)\Big\}=\sum_{j=1}^{\lfloor t/\eps\rfloor}\log(1-\alpha_{2,j}),
\]
and using \eqref{eq:assum-asym}, we get
\[
\delta_\eps(t)=\sum_{j=1}^{\lfloor t/\eps\rfloor}\log\Big(1-\eps f_2(j\eps)-\tilde{\alpha}_{2,j,\eps}\eps\Big).
\]
Due to the continuity of the function $f_2$ and to the uniform limit of $\tilde \alpha$ to zero,
\begin{equation}\label{eq:limit-inte}
\lim_{\eps\to 0}\delta_\eps(t)=-\lim_{\eps\to 0}\eps \sum_{j=1}^{\lfloor t/\eps\rfloor}f_2(j\eps)=-\int_0^t f_2(u)du.
\end{equation}
Hence for any $t\ge 0$,
\[\lim_{\eps\to 0}\mathbb{P}\left(\eps (T_{2n+1}- T_{2n})>t\right)=\exp\left(-\int_0^t f_2(u)du\right).\]
The same arguments lead to 
\[\lim_{\eps\to 0}\mathbb{P}\left(\eps (T_{2n+2}- T_{2n+1})>t\right)=\exp\left(-\int_0^t f_1(u)du\right).\] We conclude that $\mathcal{R}_n^\eps$ converges in distribution towards $(e_1,e_1+e_2,\ldots,e_1+e_2+\ldots+e_n)$, for any $n\ge 1$.
\paragraph{Step 2--- Duality and convergence of the counting process.} 
Let us define the following right-continuous counting process:
\begin{equation}\label{eq:count}
N^\eps(t)=\sup\{n\ge 0:\, \eps T_n\le t\}=\sum_{n\ge 1}\ind{\eps T_n \le t}.
\end{equation}
In order to prove \eqref{eq:skoro} we first point out the convergence of the counting process $N^\eps$ towards $N^0$. 
The one-to-one correspondence between $(N^\eps(t))_{t\geq 0}$ and $(T_n)_{n\ge 1}$ implies that for any $0<t_1<\ldots<t_k$, the convergence in distribution
of $(N^\eps(t_1),\ldots,N^\eps(t_n))$ as $\eps$ tends to zero
is a consequence of the convergence of $\rond R_n^\eps$. Indeed
\[\mathbb{P}(N^\eps(t_1)=j_1,\ldots,N^\eps(t_n)=j_n)=\mathbb{P}(\eps T_{j_1}\le t_1<\eps T_{j_1+1},\ldots, \eps T_{j_n}\le t_n<\eps T_{j_n+1})\]
and consequently
\[\lim_{\eps \to 0}\mathbb{P}(N^\eps(t_1)=j_1,\ldots,N^\eps(t_n)=j_n)=\mathbb{P}\left(E_{j_1}\le t_1<E_{j_1+1},\ldots, E_{j_n}\le t_n<E_{j_n+1}\right),\]
where $E_n=\sum_{k=1}^n e_k$. In order to obtain the convergence of the counting processes, it suffices to use a tightness criterium (see, for instance, \cite[Theorem 15.2 p. 125]{Billingsley}).
Let $s<t$ and let us denote $\tau_{st}:= \lfloor t/\eps\rfloor-\lfloor s/\eps\rfloor$ then
\begin{align}\label{eq:introd-d}
d_{s,t}^\eps &:=\mathbb{P}(N^\eps(t)>N^\eps(s))=1-\mathbb{P}(N^\eps(t)=N^\eps(s))\nonumber\\
&\ =1-\mathbb{P}(N^\eps(t)=N^\eps(s), \, N^\eps(s)\in 2\mathbb{N})-\mathbb{P}(N^\eps(t)=N^\eps(s), \, N^\eps(s)\in 2\mathbb{N}+1).
\end{align}
Since $X_0=1$, if $N^\eps(s)\in 2\mathbb{N}$ we have on one hand $X_{\lfloor s/\eps \rfloor}=1$ and on the other hand $M_{\lfloor s/\eps \rfloor}\le \lfloor s/\eps \rfloor+1$. Assuming $M_{\lfloor s/\eps \rfloor}=\ell+1$ with $0\le \ell\le \lfloor s/\eps \rfloor$ then 
\begin{align}\label{eq:prob-condi}
P_{st}(\ell)&:=\mathbb{P}\Big(N^\eps(t)=N^\eps(s)\Big|M_{\lfloor s/\eps \rfloor}=\ell+1, \, X_{\lfloor s/\eps \rfloor}=1\Big)\nonumber\\
&\ =\mathbb{P}\Big(  X_{\lfloor s/\eps \rfloor+1}=1,\ldots,X_{\lfloor s/\eps \rfloor+\tau_{st}}=1 \Big|M_{\lfloor s/\eps \rfloor}=\ell+1, \, X_{\lfloor s/\eps \rfloor}=1\Big)\nonumber\\
&\ =(1-\alpha_{2,\ell+1})(1-\alpha_{2,\ell+2})\ldots (1-\alpha_{2,\ell+\tau_{st}}).
\end{align}
Then it comes,
\begin{align}\label{eq:develop}
\mathbb{P}(N^\eps(t)=N^\eps(s), \, N^\eps(s)\in 2\mathbb{N})&=\sum_{\ell=0}^{\lfloor s/\eps \rfloor}\mathbb{P}\Big(N^\eps(t)=N^\eps(s), \, N^\eps(s)\in 2\mathbb{N},\, M_{\lfloor s/\eps \rfloor}=\ell+1\Big)\nonumber\\
&=\sum_{\ell=0}^{\lfloor s/\eps \rfloor}\prod_{k= 1}^{\tau_{st}}(1-\alpha_{2,k+\ell})\mathbb{P}\Big(N^\eps(s)\in 2\mathbb{N},\, M_{\lfloor s/\eps \rfloor}=\ell+1\Big)\nonumber\\
&\ge \inf_{0\le \ell\le \lfloor s/\eps \rfloor}\prod_{k=
1}^{\tau_{st}}(1-\alpha_{2,k+\ell})\mathbb{P}(N^\eps(s)\in
2\mathbb{N}).
\end{align}
Similar arguments are used in the odd case $N^\eps(s)\in 2\mathbb{N}+1$. In this situation $X_{\lfloor s/\eps \rfloor}=-1$ and the sequence $(\alpha_{2,\bullet})$ in \eqref{eq:prob-condi} is therefore replaced by $(\alpha_{1,\bullet})$. Combining \eqref{eq:develop} with \eqref{eq:introd-d}, we obtain
\[
d_{s,t}^\eps \le1-\inf_{0\le l\le \lfloor s/\eps \rfloor}\prod_{k= 1}^{\tau_{st}}(1-\alpha_{2,k+l})\mathbb{P}(N^\eps(s)\in 2\mathbb{N})-\inf_{0\le \ell\le \lfloor s/\eps \rfloor}\prod_{k=
1}^{\tau_{st}}(1-\alpha_{1,k+\ell})\mathbb{P}(N^\eps(s)\in
2\mathbb{N}+1).\\
\]
By \eqref{eq:assum-asym}, we get
\begin{align*}
d_{s,t}^\eps &\le 1-\inf_{i=1,2}\left\{ \inf_{0\le \ell\le \lfloor s/\eps \rfloor}\prod_{k=1}^{\tau_{st}}\Big(1-\eps f_i(\eps(k+\ell)\Big)\right\}+o(\eps)\\
&\le 1-\inf_{i=1,2}\left\{\Big(1-\eps \sup_{0\le u\le t+\eps}f_i(u) \Big)^{\tau_{st}} \right\}+o(\eps)\\
&\le 1-\Big(1-\eps \sup_{0\le u\le t+\eps} f_1(u)\vee f_2(u) \Big)^{\tau_{st}}+o(\eps) .
\end{align*}
Since $\eps \tau_{st}\le t-s+\eps$, for any $\delta>0$,
$N>0$, we can find $\eps_0>0$ such that $d_{s,t}^\eps\le
\delta$ for all $\eps\le\eps_0$ and $t,s\le N$. We deduce
that the set of all the distributions of $N^\eps$, $\eps\in]0,1]$, is weakly relatively compact and
obtain finally the convergence in law of $N^\eps$ towards
$N^0$.
\paragraph{Step 3--- \mathversion{bold}Convergence of $(S^\eps,X^\eps,M^\eps)$.\mathversion{normal}}
We have just proved that $(N^\eps(t))_{t\ge 0}$ converges in
distribution towards $(N^0(t))_{t\ge0}$. The paths of these processes belong to the Skorohod space $\mathbb{D}$. The two main ingredients of the proof are the following. First we note that $S^\eps(t)$, $X^\eps(t)$ and $M^\eps(t)$ can be expressed continuously in terms of the process $(N^\eps(s), \ s\le t)$ and secondly we use the convergence of $N^\eps$. For the process $S^\eps(t)$, we introduce the mapping $F_1:\mathbb{D}(0,1)\to\mathcal{C}(0,1)$ defined for $t\in[0,1]$ by
\[
F_1(f)(t)=\int_0^t\cos(\pi f(s))\,ds.
\]
%This mapping is continuous (see Lemma~\ref{lem:sko1} and Lemma~\ref{lem:sko2}), therefore the process $\Phi(N^\eps)$
%converges in distribution towards the process $\Phi(N^0)$.
Since $N^\eps$ is $\mathbb{N}$-valued, we get
\[
F_1(N^\eps)(t)=\int_0^t\cos(\pi N^\eps(s))\,ds=\int_0^t (-1)^{N^\eps(s)}\,ds.
\]
Note that \eqref{eq:defdeseps} combined with \eqref{eq:rel-S-N} imply that
$S^\eps(t)=\eps+\int_0^t (-1)^{N^\eps(s+\eps)} ds$.
Finally the definition of $S^\eps(t)$ leads to
\begin{equation}\label{approx1}
\vert S^\eps(t)-F_1(N^\eps)(t)\vert=\Big|\eps+\int_{0}^t(-1)^{N^\eps(s+\eps)}\,ds-\int_0^t (-1)^{N^\eps(s)}\,ds\Big|\le 3\eps.
\end{equation}
For the process $X^\eps$, we observe that 
$X^\eps(t)=F_2(N^\eps(t)):=\cos(\pi N^\eps(t))$ and the memory process is linked to the age process of $N^\eps$:
\[
\Big|M^\eps(t)-\Big(t-\inf\{s\ge 0:\ N^\eps(s)=N^\eps(t)\}\Big)\Big|\le \eps.
\]
Let us just note that for $\Phi(x)=\cos(\frac{\pi}{2} x)\ind{[-1,1]}(x)$ which is a continuous function, we get
\[
t-\inf\{s\ge 0:\ N^\eps(s)=N^\eps(t)\}=\int_0^t \Phi(N^\eps(t)-N^\eps(s))\,ds=F_3(N^\eps)(t)
\]
where
\[ F_3:f\to \Big(\int_0^t\Phi(f(t)-f(s))ds,\ t\ge 0\Big).
\]
In order to prove \eqref{eq:skoro}, it suffices to use the convergence in distribution of $N^\eps$ towards $N^0$
developed in Step 2 and the continuity in the Skorohod topology of the three
functions $F_1$, $F_2$ and $F_3$ (see Lemma \ref{lem:sko1}, \ref{lem:sko2} and \ref{lem:sko3}).
Finally we note that $(-1)^{N^0(t)}=F_1(N^0(t))$, $S^0(t)=F_2(N^0(t))$ and $m(t)=F_3(N^0(t))$.
\QED
%
%The process $(S^0(t),N^0(t),m(t))$ is Markovian and the behavior of the paths is the following:
%\begin{itemize}
%\item If $N^0(t)\in2\mathbb{N}$ and $m(t)=T$ then the slope is equal to $1$ and the next change of slope shall appear at time $t+\tau$ where $\mathbb{P}(\tau>s)=\exp-\int_0^sg(u+T)du$.
%\item If $N^0(t)\in 2\mathbb{N}+1$ and $=T>0$ then the slope is equal to $-1$ and the next change of slope shall appears at time $t+\tau$ with $\mathbb{P}(\tau>s)=\exp-\int_0^sf(u+T)du$.
%\end{itemize}
\noindent {\bf Examples.} For some particular $f_1$, the related random variable $e_{2n}$ has a distribution which belongs to well-known families of laws.
\begin{itemize}
\item If $f_1$ is a constant function then the sequence $(e_{2n})$ is exponentially distributed.
\item If $f_1(x)=\alpha\lambda x^{\alpha-1}$ with $\alpha>0$ and $\lambda>0$ then the law of $e_{2n}$ corresponds to
the Weibull distribution with parameters $(\alpha,\lambda)$.
\item If $f_1(x)=\frac{\lambda}{x}\,\ind{x\ge x_0}$ with $x_0>0$, then we deal with the Pareto
distribution for $e_{2n}$.
\end{itemize}
It has been shown in \cite{Herrmann-Vallois} that the density part of the distribution of $S(t)$ can be expressed via Bessel functions. Here, we have a weaken result which says that we are only able to determine the Laplace transform of $S(\tau)$ (see, Proposition \ref{prop:double-lapl} below). Being unable to invert this transformation, the distribution of $S(t)$ is unknown. Although the path description of $(S(t))_{t\ge 0}$ is very easy, only few properties related to the GITN are known.
%
%
%
%
%The description of the limit process is of particular interest. An explicit expression based on the Bessel functions can be pointed out for the classical integrated telegraph noise \cite{Herrmann-Vallois}. For the generalized situation, similar computations can not be carried through since we do not have the same Markovian properties. Nevertheless the convergence of the persistent random walk towards the GITN permits to use the generating function described in Theorem \ref{thm:doublelap} to obtain an explicit expression
%for the double Laplace transform of $S^0(t)$ for $t\ge 0$.}
\begin{prop}\label{prop:double-lapl}
Let $(S^0(t))_{t\in \mathbb{R}_+}$ be the GITN defined by \eqref{def:gitn} then the double Laplace transform
defined by
\begin{equation}\label{eq:dbl}
\mathcal{L}(r,\gamma):=\int_0^\infty e^{-rt}\ \mathbb{E}\Big[ e^{-\gamma S^0(t)} \Big]\,dt,\quad r>0,\ \gamma>0,
\end{equation}
is equal to
\[
\frac{-(r+\gamma)\mathcal{R}(r-\gamma,f_1)\mathcal{R}(r+\gamma,f_2)+\mathcal{R}(r-\gamma,f_1)+\mathcal{R}(r+\gamma,f_2)}{(r-\gamma)\mathcal{R}(r-\gamma,f_1)+(r+\gamma)\mathcal{R}(r+\gamma,f_2)-(r^2-\gamma^2)\mathcal{R}(r-\gamma,f_1)\mathcal{R}(r+\gamma,f_2)},
\]
where
\begin{equation}\label{eq:dbl2}
\mathcal{R}(z,f_i)=\int_0^\infty e^{-zt-\int_{0}^t f_i(u)\,du}\,dt,\quad z\in\mathbb{R},\ i=1,2.
\end{equation}
\end{prop}
\begin{rem}\begin{enumerate}\item[(i)]
In the particular constant case, that is $f_1(t)=f_1$ and $f_2(t)=f_2$ for all $t\ge 0$, the stochastic process
corresponds to the so-called integrated telegraph noise introduced in \cite{Herrmann-Vallois}. For this process,
we get
\(
\mathcal{R}(z,f_i)=(z+f_i)^{-1}\) for  $i=1,2$.
The double Laplace transform $\mathcal{L}$ becomes
\[
\mathcal{L}(r,\gamma)=\frac{f_0+g_0+r-\gamma}{r^2-\gamma^2+(r-\gamma) g_0+(r+\gamma)f_0}.
\]
This identity was already obtained by Weiss in \cite{weiss02} and presented in \cite{Herrmann-Vallois} (see Remark 3.10).
\item[(ii)] Let $\xi$ be an exponential r.v. with parameter $r$ independent from $(S^0(t),\ t\ge 0)$. Then $\mathcal{L}(r,\gamma)$ is the Laplace transform of $S^0(\xi)$: \[\mathcal{L}(r,\gamma)=\mathbb{E}[e^{-\gamma S^0(\xi)}].
\]
\end{enumerate}
\end{rem}
\noindent {\sc Proof of Proposition \ref{prop:double-lapl}.}
Recall that $S^\eps(t)$ is the piecewise continuous process defined by \eqref{eq:defdeseps}. By Theorem \ref{thm:convergenceGITN} and the Lebesgue convergence theorem,
we just need to study the convergence of $\mathcal{L}^\eps(r,\gamma)$ the double Laplace transform of $S^\eps(t)$. As $\eps\to 0$, we get
\begin{align}\label{eq:approxlapl}
\mathcal{L}^\eps(r,\gamma)&=\int_0^\infty e^{-rt}\ \mathbb{E}\Big[ e^{-\gamma S^\eps(t)} \Big]\,dt=\sum_{k\ge 0}\int_{k\eps}^{(k+1)\eps} e^{-rt}\Big( \mathbb{E}[e^{-\gamma \eps S_k}]+o(\eps) \Big)\,dt\nonumber\\
&= \frac{1-e^{r\eps}}{r}\sum_{k\ge 0}\Big(e^{-r\eps}\Big)^k\mathbb{E}\Big[ e^{-\gamma\eps S_k} \Big]+o(\eps)=\frac{1}{r} \mathbb{E}[(e^{-\gamma \eps})^{S_\tau}]+o(\eps),
\end{align}
where $\tau+1$ is a geometrically distributed random variable, independent of the process $(S_n)$:
\[
\mathbb{P}(\tau=n)=(e^{-r\eps})^n(1-e^{-r\eps}).
\]
Obviously \eqref{eq:approxlapl} shows that $r\mathcal{L}^\eps(r,\gamma)$ and $\mathbb{E}[e^{-\gamma \eps S_\tau}]$ have the same limit as $\eps \to 0$.
Note that choosing $\lambda=e^{-\gamma \eps}$ and $\rho=e^{-r\eps}$ in Theorem \ref{thm:doublelap} gives the value of $\mathbb{E}[e^{-\gamma\eps S_\tau}]$. Due to the specific form of \eqref{eq:thm:double} we are lead to prove the following intermediate result:
\begin{align}\label{eq:fin}
\lim_{\eps\to 0}\eps\widehat{\mathcal{P}}_i(e^{-\eps z})=\mathcal{R}(z,f_i).
\end{align}
where $\hat{\mathcal{P}}_i$ (resp. $\mathcal{R}(z,f_i)$) is defined by \eqref{eq:def:pchap} (resp. \eqref{eq:dbl2}). \\
Indeed, according to the definition of $\hat{\mathcal{P}}_i$ we easily get
\[
\eps\widehat{\mathcal{P}}_i(e^{-\eps z})=\frac{z\eps}{1-e^{-\eps z}}\ e^{-\eps z}\int_0^\infty \mathcal{P}_i\Big( \left\lfloor \frac{t}{\eps} \right\rfloor +1 \Big)e^{-zt}dt.
\]
Using \eqref{eq:limit-inte} (where the index 2 is replaced by $i$) yields 
\begin{align*}
\lim_{\eps\to 0}\mathcal{P}_i\Big( \left\lfloor \frac{t}{\eps} \right\rfloor +1 \Big)=\lim_{\eps\to 0} e^{\delta_\eps(t)}=e^{-\int_{0}^t f_i(u)\,du}.
\end{align*}
Then, the dominated convergence theorem implies \eqref{eq:fin}. Since
\begin{itemize}
\item $S^\eps(t)$ converges in distribution to $S^0(t)$ as $\eps\to 0$
\item $\rho-1\sim -r\eps$ and $\lambda \rho -1\sim -(r+\gamma)\eps$ as $\eps \to 0$
\end{itemize}
then \eqref{eq:approxlapl} and Theorem \ref{thm:doublelap} imply
\begin{align*}
&\mathcal{L}(r,\gamma)=\lim_{\eps\to 0}\mathcal{L}^\eps(r,\gamma)=\frac{1}{r}\lim_{\eps\to 0}\mathbb{E}[(e^{-\gamma \eps})^{S_\tau}]\\
&\ =\lim_{\eps\to 0} \frac{-(r+\gamma)R_1^\eps(r-\gamma)R_2^\eps(r+\gamma)+R_1^\eps(r-\gamma)+R_2^\eps(r+\gamma)}{(r-\gamma)R_1^\eps(r-\gamma)+(r+\gamma)R_2^\eps(r+\gamma)-(r^2-\gamma^2)R_1^\eps(r-\gamma)R_2^\eps(r+\gamma) }
\end{align*}
where $R_1^\eps(z)=\eps\widehat{\mathcal{P}}_1(e^{-z\eps})$.\\
It is clear that Proposition \ref{prop:double-lapl} is a straightforward consequence of \eqref{eq:fin} and the above identity.
\QED
\appendix
\section{Continuity in the Skorohod space}
\setcounter{equation}{0}
Let us denote $\mathbb{D}([0,1])$ the Skorohod space \emph{i.e.} the space of functions which are
right-continuous and have left-hand limits. $\mathbb{D}$ is a complete metric space for the following distance
(see \cite[Theorem 14.2]{Billingsley})
\begin{equation}\label{eq:dist-sko}
d(f,g)=\inf_{\lambda\in\Lambda}\max\Big\{ \Vert \lambda \Vert,\, \Vert f-g\circ \lambda\Vert_\infty \Big\},
\end{equation}
where
\[
 \Vert \lambda \Vert=\sup_{s\neq t}\left|\log\frac{\lambda(t)-\lambda(s)}{t-s}\right|,
\]
$\Vert \cdot\Vert_\infty$ is the uniform norm and $\Lambda$ is the space of strictly increasing,
continuous mappings of $[0,1]$ into itself.
\begin{lem}\label{lem:sko1}
Let $\Phi:\mathbb{R}\to\mathbb{R}$ be a continuous function, then $f\in \mathbb{D}([0,1])\to \Phi \circ f$ is continuous in the Skorohod topology.
\end{lem}
\pff
Let $f\in  \mathbb{D}([0,1])$. Then there exists $M>0$ such that $|f(t)|\le M$ for all $t\in[0,1]$. For $\eps>0$,
due to the uniform continuity of $\Phi$, there exists $\delta>0$ such that: for any $(x,y)\in[-2M,2M]^2$
satisfying $|x-y|<\delta$ we have
$|\Phi(x)-\Phi(y)|<\eps$.
Let us consider now a function $g\in  \mathbb{D}([0,1])$ such that $d(f,g)<\delta\wedge M$. Therefore, there exists $\lambda\in\Lambda$ such that $\Vert \lambda\Vert_\infty<\delta$ and $\Vert f-g\circ\lambda\Vert_\infty<\delta$. Consequently
\[
\Vert \Phi(f)-\Phi(g\circ\lambda)\Vert_\infty< \eps.
\]
Continuity of $\Phi$ at $f$ follows from the definition of Skorohod's distance.
\QED
\begin{lem}\label{lem:sko2}
The mapping $f\in  \mathbb{D}([0,1])\to \Big(\int_0^t f(u)\,du,\ t\ge 0\Big)$ is continuous in the Skorohod topology.
\end{lem}
\pff First let us recall that any function belonging to the Skorohod space is integrable.
We denote $I_f(t)=\int_0^t f(u)\,du$. Let $f,g\in  \mathbb{D}([0,1])$ such that $d(f,g)<\delta$ and choose $\lambda\in\Lambda$
with $\Vert \lambda\Vert<\delta$ and $\Vert f-g\circ \lambda\Vert_\infty<\delta$, we get
\begin{align}\label{eq:skoint}
\vert I_f(t)-  I_g\circ\lambda(t) \vert&=\lim_{n\to\infty}\left|\sum_{k=1}^n \frac{t}{n} f\Big(\frac{kt}{n}\Big) -g\circ\lambda\Big(\frac{kt}{n}\Big)\Big\{ \lambda\Big(\frac{kt}{n}\Big)-\lambda\Big(\frac{(k-1)t}{n}\Big) \Big\}\right|\nonumber \\
&\le \lim_{n\to\infty}\left|\frac{t}{n}\sum_{k=1}^n (f-g\circ\lambda)\Big(\frac{kt}{n}\Big)\right|.\nonumber \\
&+\lim_{n\to\infty}\left|\sum_{k=1}^ng\circ\lambda\Big(\frac{kt}{n}\Big)  \Big\{ \lambda\Big(\frac{kt}{n}\Big)-\lambda\Big(\frac{(k-1)t}{n}\Big) -\frac{t}{n}\Big\} \right|.
\end{align}
By definition of the norm on the Skorohod space, we have
\[
e^{-\Vert \lambda \Vert}<\frac{\lambda(t)-\lambda(s)}{t-s}<e^{\Vert \lambda \Vert}, \quad\mbox{for}\ 0\le s<t\le 1.
\]
Consequently for any $0\le s<t\le 1$, we have
\begin{equation}
\label{eq:bound-norm}
|\lambda(t)-\lambda(s)-(t-s) |\le (t-s)\max\Big(e^{\Vert \lambda\Vert}-1,1-e^{-\Vert\lambda\Vert}\Big)\le (t-s)(e^{\Vert \lambda\Vert}-1).
\end{equation}
Combining \eqref{eq:skoint} and \eqref{eq:bound-norm} yields to
\[
\vert I_f(t)-  I_g\circ\lambda(t) \vert\le \Vert f-g\circ \lambda\Vert_\infty+ \Vert g\Vert_\infty (e^{\Vert \lambda\Vert} -1),\quad \forall t\in [0,1].
\]
We deduce that $d(f,g)<\delta$ implies
\[
d(I_f,I_g)\le \max\{ \delta, \delta+ \delta e^\delta \Vert g\Vert_\infty \}= \delta (1+e^\delta \Vert g\Vert_\infty).
\]
As a result $f\to \int_0^1 f(u)du$ is a continuous mapping.
\QED
Using similar arguments as those presented in the proofs of Lemma~\ref{lem:sko1} and Lemma~\ref{lem:sko2},
we obtain the following continuity result.
\begin{lem}\label{lem:sko3}
Let $\Phi$ be a continuous function, then the mapping
\[
f\in  \mathbb{D}([0,1])\longrightarrow \left( \int_0^t \Phi(f(t)-f(s))\,ds,\ t\ge 0 \right)
\]
is continuous in Skorohod's topology.
\end{lem}
%%%%%%%%%%%%%%%%%%%%%%%%%%%%%
\section{Invariant measure for the double infinite comb}\label{app:double-peigne}
\setcounter{equation}{0}
Consider the probabilized context tree given on Figure~\ref{figPeigneInfini}. In this case, there are two infinite leaves $0^\infty$ and $1^{\infty}$ and a countable number of leaves $0^n1$ and $1^n0$, $n\in\mathbb{N}$. Suppose that $\pi$ is a stationary measure on $\mathcal{L}$. Denote by $\rond W$ the set of finite words on the alphabet $\{ 0,1\}$. For any finite word $w\in \rond W$, we denote by $\pi (w): = \pi(\rond Lw)$ the measure of the cylinder $\rond Lw$ denoting the set of left infinite words ending with $w$.
We first compute $\pi (w)$ as a function of $\pi (1)$ when the reversed word of $w$ is any context or any internal node. Applying equation~(\ref{eq:def:VLMC}) to $U_n=\ldots 10^n$, it comes for any $n\geq 1$,
\[\pi(10^n)=\pi(10^{n-1})q_{0^{n-1}1}(0).\]
An immediate induction yields, for any $n\geq 1$,
\begin{equation}
\label{probaContextesDoublePeigneGauche}
\pi (10^n)=\pi (10)\prod_{k=1}^{n-1}q_{0^k1}(0)=\pi (10)\prod_{k=1}^{n-1}(1-\alpha_{1,k})=\pi (10)\cP_1(n).
\end{equation}
In the same way, 
\begin{equation}
\label{probaContextesDoublePeigneDroite}
\pi (01^n)=\pi(01)\cP_{2}(n),
\end{equation}
The stationary probability of a reversed context is thus necessarily given by
Formulae~(\ref{probaContextesDoublePeigneGauche}) and (\ref{probaContextesDoublePeigneDroite}).
Now, if $0^n$ is any internal node of the context tree but $0$, we need going down along the branch in the context tree to reach the contexts; using then the disjoint union $\pi (0^{n+1})=\pi (0^n)-\pi (10^n)$, by induction, it comes for any $n\geq 2$, 
\begin{equation}
\label{probaNoeudsInternesGauchesDoublePeigne}
\pi (0^n)=\pi(0)-\pi(10)\sum _{k=1}^{n-1}\cP_{1}(k).
\end{equation}
The same holds for any internal node $1^n$ but $1$,
\begin{equation}
\label{probaNoeudsInternesDroitesDoublePeigne}
\pi (1^n)=\pi(1)-\pi(10)\sum _{k=1}^{n-1}\cP_{2}(k).
\end{equation}
where we have used $\pi(01)=\pi(10)$ (coming from the invariance of $\pi$).
The stationary probability of a reversed internal node of the context tree is thus necessarily given by
Formulae~(\ref{probaNoeudsInternesGauchesDoublePeigne}) and (\ref{probaNoeudsInternesDroitesDoublePeigne}).

\vskip 5pt
It remains to compute $\pi (10)$ and then $\pi(0)$ (and consequently $\pi(1)$).
The denumerable partition of the whole probability space given by all cylinders based on leaves in the context
tree implies
$1-\pi (0^\infty )-\pi(1^\infty)=\pi (10)+\pi (100)+\dots+\pi(01)+\pi(011)+\ldots$, \ie
\begin{equation}
\label{partitionDoublePeigne}
1-\pi (0^\infty )-\pi(1^\infty)=\pi (10)\sum_{n\geq 1}\left(\cP_{1}(n)+\cP_{2}(n)\right).
\end{equation}

This leads to the following statement that covers all cases of existence, unicity and nontriviality for a stationary
probability measure for the double infinite comb.
In the generic case (named \emph{irreducible} case hereunder), we give a necessary and sufficient condition on the data for the existence of a stationary probability measure;
moreover, when a stationary probability exists, it is unique.
The \emph{reducible} case is much more singular and gives rise to nonunicity. 

\begin{prop}
\label{mesureStationnaireDoublePeigne}
{\bf (Stationary probability measures for a double infinite comb)}

Let $(U_n)_{n\geq 0}$ be a VLMC defined by a probabilized double infinite comb.

\medskip
\begin{enumerate}[(i)]
\item \emph{Irreducible case}: Assume that $q_{0^\infty}(0)\neq 1$ and $q_{1^\infty}(1)\neq 1$.
\begin{enumerate}[(a)]
\item \emph{Existence}: The Markov process $(U_n)_{n\geq 0}$ admits a stationary probability measure on $\mathcal{L}$ if and only if
the numerical series $\Theta_1$ and $\Theta_2$ converge.
\item \emph{Unicity}: Assume that the series $\Theta_1$ and $\Theta_2$ converge. Then, the stationary probability measure $\pi$ on $\mathcal{L}$ is unique;
it is characterized by 
\begin{equation}
\label{pi(1)DoublePeigne}
\pi (0)=\frac{\Theta_1}{\Theta_1+ \Theta_2}\ \ ,\ \ 
\pi(10)=\frac{1}{\Theta_1+\Theta_2}
\end{equation}
and Formulae~(\ref{probaContextesDoublePeigneGauche}), (\ref{probaContextesDoublePeigneDroite})
(\ref{probaNoeudsInternesGauchesDoublePeigne}), (\ref{probaNoeudsInternesDroitesDoublePeigne}).
\end{enumerate}
\item \emph{Reducible cases}: Assume that $q_{0^\infty}(0)=1$ and $q_{1^\infty}(1)\neq 1$.
\begin{enumerate}[(a)]
\item If at least one of the series $\Theta_1$ and $\Theta_2$ diverges, then the trivial probability measure $\pi$ on $\mathcal{L}$ defined by
$\pi (0^\infty )=1$ is the unique stationary probability measure.
\item If the series $\Theta_1$ and $\Theta_2$ converge, then there is a one parameter family of stationary probability measures
on $\mathcal{L}$.
More precisely, for any $a\in [0,1]$, there exists a unique stationary probability measure $\pi _a$ on $\mathcal{L}$
such that $\pi _a(0^\infty)=a$.
The probability $\pi _a$ is characterized by
\[\pi _a(0)=\frac {a\Theta_2+\Theta_1}{\Theta_1+\Theta_2}, \quad \pi_a(10)=\frac{1-a}{\Theta_1+\Theta_2}\] 
and Formulae~(\ref{probaContextesDoublePeigneGauche}), (\ref{probaContextesDoublePeigneDroite})
(\ref{probaNoeudsInternesGauchesDoublePeigne}), (\ref{probaNoeudsInternesDroitesDoublePeigne}). 
\smallskip

Assume that $q_{0^\infty}(0)\neq 1$ and $q_{1^\infty}(1)=1$. Then the same results as in (ii.a) and (ii.b) hold, exchanging the role of $0$ and $1$.

\smallskip

Assume that $q_{0^\infty}(0)=1$ and $q_{1^\infty}(1)=1$. 

\item If at least one of the series $\Theta_1$  and $\Theta_2$ diverges, then there is a one parameter family of stationary probability measures
on $\mathcal{L}$.
More precisely, for any $a\in [0,1]$, there exists a unique stationary probability measure $\pi _a$ on $\mathcal{L}$
such that $\pi _a(0^\infty)=a$.
The probability $\pi _a$ is characterized by
$\pi_a(0^n)=a$ and $\pi_a(1^n)=1-a$ for every $n\geq 1$ and $\pi _a(w)=0$ as soon as $w$ contains one $0$ and one $1$.
\item If the series $\Theta_1$ and $\Theta_2$ converge, then there is a two parameters family of stationary probability measures
on $\mathcal{L}$.
More precisely, for any $a\in [0,1]$ and $b\in [0,1]$, there exists a unique stationary probability measure $\pi _{a,b}$ on $\mathcal{L}$
such that $\pi _{a,b}(0^\infty)=a$ and $\pi _{a,b}(1^\infty)=b$.
The probability $\pi _{a,b}$ is characterized by
\[\pi _{a,b}(0)=\frac {a\Theta_2+(1-b)\Theta_1}{\Theta_1+\Theta_1},\quad \pi_{a,b}(10)=\frac{1-a-b}{\Theta_1+\Theta_2}\]
and Formulae~(\ref{probaContextesDoublePeigneGauche}), (\ref{probaContextesDoublePeigneDroite})
(\ref{probaNoeudsInternesGauchesDoublePeigne}), (\ref{probaNoeudsInternesDroitesDoublePeigne}). 
\end{enumerate}
\end{enumerate}
\end{prop}

\pff
%\begin{enumerate}[\emph{(i)}]
\begin{enumerate}[(i)]
\item Assume that $q_{0^\infty}(0)\neq 1$, $q_{1^\infty}(1)\neq 1$ and that $\pi$ is a stationary probability measure.
By definition of probability transitions, $\pi (0^\infty )=\pi (0^\infty )q_{0^\infty}(0)$ and $\pi (1^\infty )=\pi (1^\infty )q_{1^\infty}(1)$ so that $\pi (0^\infty )$
and $\pi(1^\infty)$ necessarily vanish.
Thus, thanks to~(\ref{partitionDoublePeigne}), $\pi (10)\neq 0$, the series $\Theta_1+\Theta_2$ converges and so do $\Theta_1$ and $\Theta_2$. This also implies 
\[1=\pi(10)(\Theta_1+\Theta_2).\]
Passing to the limit in~(\ref{probaNoeudsInternesGauchesDoublePeigne}) implies $\pi(0)=\pi(10)\Theta_1$.
Thus Formula~(\ref{pi(1)DoublePeigne}) is valid.
Moreover, when $\overline w$ is any context or any internal node of the context tree, $\pi (w)$ is
necessarily given by Formulae~(\ref{pi(1)DoublePeigne}), (\ref{probaContextesDoublePeigneGauche}), (\ref{probaContextesDoublePeigneDroite}), (\ref{probaNoeudsInternesGauchesDoublePeigne})
and~(\ref{probaNoeudsInternesDroitesDoublePeigne}). Since the cylinders $\rond L w$, $w\in\rond W$ span the $\sigma$-algebra on $\rond L$, there is at most one stationary
probability measure.
This proves the \emph{only if} part of {\it (i.a)}, the unicity and the characterization claimed in~{\it (i.b)}.

Reciprocally, when the series converge, Formulae~(\ref{pi(1)DoublePeigne}), (\ref{probaContextesDoublePeigneGauche}), (\ref{probaContextesDoublePeigneDroite})
(\ref{probaNoeudsInternesGauchesDoublePeigne}), (\ref{probaNoeudsInternesDroitesDoublePeigne}) define a probability measure on
the semiring spanned by cylinders, which extends to a stationary probability measure on the whole
$\sigma$-algebra on $\rond L$.
This proves the \emph{if} part of {\it (i.a)}.

To deal with the reducible cases, recall the three following equations (which hold when the series converge) :
\[\left\{\begin{array}{rcl}\label{systeme}
  1-\pi(0^\infty)-\pi(1^\infty) & = & \pi(10)(\Theta_1+\Theta_2) \\
  \pi(0^\infty) & = & \pi(0)-\pi(10)\Theta_1 \\
  \pi(1^\infty) & = & \pi(1)-\pi(10)\Theta_2
  \end{array}\right.\]
\item 
Assume that $q_{0^\infty}(0)=1$ and $q_{1^\infty}(1)\neq 1$ .
First, as above, $q_{1^\infty}(1)\neq 1$ implies $\pi(1^\infty)=0$.
Next, Formula~(\ref{partitionDoublePeigne}) is always valid so that the divergence of at least one of the series forces $\pi (10)$ to vanish. This gives $\pi(0^\infty)=1$.
With the assumption $q_{0^\infty}(0)=1$, one immediately sees that this trivial probability is stationary,
proving {\it (ii.a)}.

To prove {\it (ii.b)}, assume furthermore that the series $\Theta_1$ and $\Theta_2$ converge and let $a\in [0,1]$.
As before, any stationary probability measure $\pi$ is completely determined by $\pi(0)$ and $\pi(10)$. As above, $\pi(1^\infty)=0$ and if we fix $\pi(0^\infty)=a$, the system (\ref{systeme}) reduces to
\[\left\{\begin{array}{rcl}
  1-a & = & \pi(10)(\Theta_1+\Theta_2) \\
  a & = & \pi(0)-\pi(10)\Theta_1
  \end{array}\right.\]
This gives the characterisation of {\it (ii.b)}.
Formulae~(\ref{probaContextesDoublePeigneGauche}), (\ref{probaContextesDoublePeigneDroite})
(\ref{probaNoeudsInternesGauchesDoublePeigne}), (\ref{probaNoeudsInternesDroitesDoublePeigne}) standardly extend $\pi_a$ to the
whole $\sigma$-algebra on $\rond L$ and $\pi_a$ is clearly stationary.

{\it (ii.c)}
Assume that $q_{0^\infty}(0)=1$ and $q_{1^\infty}(1)=1$ . As previously, Formula~(\ref{partitionDoublePeigne}) is valid so that the divergence of at least one of the series forces $\pi (10)$ to vanish. Let $a\in [0,1]$ and fix 
$\pi(0^\infty)=a$, the system (\ref{systeme}) reduces to $\pi(0^\infty)=\pi(0)=a$ and $\pi(1^\infty)=\pi(1)=1-a$. The invariance of this measure may be easily checked.

To prove {\it (ii.d)}, assume furthermore that the series $\Theta_1$ and $\Theta_2$ converge and let $a\in [0,1]$ and $b\in [0,1]$. If we fix $\pi(0^\infty)=a$ and $\pi(1^\infty)=b$, the system (\ref{systeme}) is equivalent to
  $$\left\{\begin{array}{rcl}
  \pi(0)-\pi(10)\Theta_1 & = & a \\
  \pi(0)+\pi(10)\Theta_2 & = & 1-b
  \end{array}\right.$$
As $\Theta_1\geq 1$ and $\Theta_2\geq 1$, this system has a unique solution given by  
\[\pi _{a,b}(0)=\frac {a\Theta_2+(1-b)\Theta_1}{\Theta_1+\Theta_2}\quad \mbox{and} \quad \pi_{a,b}(10)=\frac{1-a-b}{\Theta_1+\Theta_2}.\]
\end{enumerate}
\QED
\subsection*{Acknowledgements}

We are very grateful to F.~Paccaut and N.~Pouyanne who made the present of the calculation of the invariant measure for the double comb to us.  

\begin{small}
\bibliographystyle{plain}
\bibliography{biblio}
\end{small}
\end{document}